\documentclass[trsc,nonblindrev]{informs4}
\usepackage{eqndefns-left} %
\RequirePackage{tgtermes}
\RequirePackage{newtxtext}
\RequirePackage{newtxmath}
\RequirePackage{bm}

\OneAndAHalfSpacedXII %

\usepackage{natbib}
 \bibpunct[, ]{(}{)}{,}{a}{}{,}%
 \def\bibfont{\small}%

\EquationsNumberedThrough %

\TheoremsNumberedThrough %
\ECRepeatTheorems %

\MANUSCRIPTNO{}

\usepackage{algorithm}
\usepackage[noend]{algorithmic}

\usepackage{booktabs}
\usepackage{placeins}

\usepackage{xcolor}

\newcommand{\NGset}{\textit{NG}}
\newcommand{\DP}{\textit{DP}}
\newcommand{\Pext}{\ensuremath P^{\text{ext}}}

\begin{document}

\RUNAUTHOR{Yarkony et al.}

\RUNTITLE{Strengthening Compact Vehicle Routing Formulations}

\TITLE{Strengthening Compact Vehicle Routing Formulations through Arc-Flow Aggregation and Refinement}

\ARTICLEAUTHORS{%
\AUTHOR{Julian Yarkony\textsuperscript{a},
Willem-Jan van Hoeve\textsuperscript{b},
Anthony Karahalios\textsuperscript{c},
Martin Karyamsetty\textsuperscript{a},
Soheyl Khalilpourazari\textsuperscript{a},
Louis-Martin Rousseau\textsuperscript{d},
Xiaohang Zhu\textsuperscript{a}
}
\AFF{\textsuperscript{a}Decomposition Lab, Optym;
\textsuperscript{b}Tepper School of Business, Carnegie Mellon University;
\textsuperscript{c}Naval Postgraduate School;
\textsuperscript{d}Polytechnique Montr\'{e}al}
} %

\ABSTRACT{%
Arc-flow formulations can provide strong relaxations for mixed-integer program formulations of vehicle routing problems, but the associated networks may be prohibitively large. We study the setting in which multiple arc flow relaxations may be applied, and introduce a dual-guided aggregation and refinement framework that constructs substantially smaller arc-flow formulations while recovering the LP bound of the original networks. Dual information is used both to identify states that can be aggregated without weakening the bound and to guide node disaggregation. We show that the resulting procedure converges finitely to the arc-flow bound of the full arc-flow formulation. The framework retains the original decision variables, allowing cutting planes to be added directly in the original variable space. Moreover, it produces an explicit MILP formulation without relying on column generation or pricing to determine the arc-flow relaxation bound. We apply our approach to the vehicle routing problem with time windows for which we define two arc flow relaxations, as well as two families of cutting planes over the original variable space. Computational experiments show that strong bounds can be obtained using substantially aggregated networks. The resulting root bounds are competitive with those obtained by the branch-and-price solver RouteOpt, while retaining a formulation that can be solved directly by a general-purpose MILP solver.
}%


\maketitle

\section{Introduction}
\label{sec:introduction}

Strong discrete relaxations of mixed-integer linear programs (MILPs) for vehicle routing problems can be obtained from state-transition networks arising from dynamic programming. The associated arc-flow formulations provide extended formulations of the convex hulls of these relaxations and can substantially strengthen the original compact formulation. Their practical use, however, may be limited by the size of the underlying networks. A detailed state description produces a strong formulation, but may require a prohibitively large number of nodes and arcs.

One approach to this difficulty is to avoid representing the full discrete solution space explicitly, and apply a Dantzig-Wolfe reformulation that is solved by column generation. In this case, the pricing problem to generate new variables can be based on an implicit labeling algorithm over the dynamic program~\citep{DanWol60}. For many problems, in particular vehicle routing applications, the resulting dual bounds are often substantially tighter than the linear programming relaxation based on the standard compact formulation~\citep{Desrochers1992}. For solving the original MILP, however, one usually resorts to branch-and-price \citep{barnprice}, or branch-cut-and-price \citep{gschwind2019stabilized}. These methods provide powerful solution frameworks, with specialized components for pricing, branching, and cutting planes. For instance, branching rules must remain compatible with pricing, often requiring specialized schemes such as successor-based branching \citep{ryan1981integer}; and incorporating strong inequalities can be difficult when they do not admit a clean representation in the pricing subproblem \citep{jepsen2008subset,archetti2011column}.

On the other hand, explicit arc-flow formulations can be added directly to the compact MILP formulation together with constraints that link the arc-flow variables to the original variables.~This can work well for suitably small networks, e.g., as the state-transition graph of a dynamic program \citep{eppen1987,martin1990}, multi-valued decision diagrams \citep{Cir_et_al_2019}, and relaxed decision diagrams \citep{BerCir2018,vanhoeve2022gcdd}. For many large-scale problems, however, an explicit embedding of arc-flow models is not practical due to their large size.
 In this work, we introduce Dual-Guided Aggregation and Refinement (DGAR), an iterative aggregation and refinement procedure that yields an arc-flow formulation of small size while recovering the relaxation bound of the original arc-flow formulation.

Our contributions are as follows. First, we establish a dual condition for bound-preserving aggregation of arc-flow formulations. Second, we develop a complementary dual-guided refinement procedure that disaggregates nodes to strengthen the formulation. Third, we combine these operations into an iterative procedure and prove finite convergence to the bound defined by the original arc-flow representation. Lastly, we evaluate our approach on the vehicle routing problem with time windows, combining two discrete route relaxations with cutting planes in the original routing variables and reduced-cost arc removal. The resulting explicit formulation can be passed directly to a general-purpose MILP solver. Our computational study examines bound quality, the contribution of cutting planes, root-node effort, and formulation size, followed by the use of the strengthened formulation for integer optimization.

The remainder of the paper is organized as follows. Section~\ref{sec:related} reviews related work. Section~\ref{sec:problem} defines the mixed-integer programming formulation and arc-flow relaxation. Sections~\ref{sec:aggregation} and~\ref{sec:refinement} develop dual-guided aggregation and refinement, respectively, and Section~\ref{sec:algorithm} combines them into the overall algorithm and establishes finite convergence.
Section~\ref{sec:vrptw} applies the framework to the VRPTW, and Section~\ref{sec:experiments} reports the computational results. Section~\ref{sec:conclusion} concludes.

\section{Related Work}
\label{sec:related}

Our work is most closely related to three approaches for obtaining strong formulations:
Dantzig-Wolfe reformulations and column generation, arc-flow formulations,
and methods based on state aggregation and iterative refinement.

Dantzig-Wolfe reformulations strengthen the linear programming relaxation of an MILP by replacing the relaxation associated with a subset of its constraints by the convex hull of the corresponding discrete feasible points \citep{DanWol60}. The resulting formulation can provide a substantially stronger bound, but generally requires exponentially many variables. Column generation avoids enumerating these variables by maintaining a restricted master problem and generating improving variables, or columns, through a pricing problem. When the discrete set admits a dynamic programming representation, the pricing problem can be solved over its state-transition network, with each source-sink path representing a discrete solution \citep{cuttingstock,Desrochers1992}. For integer optimization, column generation is typically embedded within branch-and-price, and often branch-cut-and-price \citep{barnprice}. This requires branching rules compatible with pricing and may require specialized procedures for stabilization, cutting planes, and primal heuristics. This often leads to increasingly sophisticated computational frameworks, e.g., for solving vehicle routing applications \citep{pessoa2020generic}. See \citep{desrosiers2026branch} for a comprehensive treatment of branch-and-price.

Arc-flow formulations provide an alternative representation of discrete sets in which solutions are represented by source-to-sink paths in a directed network. A particularly important construction arises from dynamic programming, where the state-transition graph of a dynamic program defines the network \citep{wolsey1973,martin1990}. Arc-flow representations can also be obtained from other graphical representations, including decision diagrams \citep{BerCir2018,vanhoeve2022gcdd}, or through problem-specific network constructions; see \citet{delima2022arcflow} for a survey.
Arc-flow and Dantzig-Wolfe reformulations are closely related and provide the same relaxation bound when defined over the same discrete set \citep{yamin2026}. Arc-flow formulations replace the potentially exponential number of path variables by an explicit network, making it possible in some cases to solve the resulting formulation directly with a general-purpose MILP solver. Their effectiveness, however, depends on the size of the network, and detailed dynamic programming representations can themselves produce prohibitively large networks.

This trade-off has motivated several approaches for reducing the computational burden of arc-flow formulations. One direction is to use decomposition methods to solve the arc-flow MILPs. \citet{HosHoe2021} apply a Benders decomposition to solve a combined arc-flow formulation for the multiple sequence alignment problem as an MILP. \citet{kowalczyk_et_al_2024} use branch-and-price to solve an arc-flow model based on a decision diagrams representation for parallel machine scheduling problems. Likewise,
\citet{delima_et_al_2023} develop a computational framework for solving arc-flow models that combines column-and-row generation, reduced-cost variable fixing, and a specialized branching scheme. Our approach instead aims to minimize the size of the arc-flow formulation and leverage off-the-shelf MILP solvers.

Many works have applied state aggregations in some form to solve smaller-sized problems more efficiently. In the context of dynamic programming, \citet{bean1987aggregation} develop a state aggregation and disaggregation procedure find solutions with a guaranteed optimality error bound. For pricing problems in the context of column generation for vehicle routing problems, state-space relaxations have been instrumental in providing strong bounds efficiently \citep{christofides1981state,baldacci2011new}; these relaxations define aggregate states implicitly by re-defining the dynamic program. Relaxed decision diagrams, on the other hand, apply a node merging function that aggregates states and which yield dual bounds for discrete optimization problems~\citep{bergman2016}. In a related network-reduction setting, \citet{BerBodCarCir2022} construct dynamic programming networks for multiobjective discrete optimization and develop node-merging and arc-removal operations that reduce these networks while preserving the Pareto frontier.

Lastly, iterative refinement has been used for constructing network representations of smaller size than the exact state-space.
Decremental state-space relaxations improve column generation by incrementally strengthening a relaxed pricing problem \citep{righini2008new}.
Dynamic discretization discovery iteratively refines a state-space network by discretizing state values such as time~\citep{boland2017continuous}.
\citet{clautiaux2017iterative} present an iterative aggregation and disaggregation method to solve arc-flow formulations by strengthening an initial relaxed arc-flow formulation.
Lastly, column elimination iteratively removes infeasible paths from a decision diagram and more generally computes locally refined dynamic programming state-transition graphs as arc-flow formulations \citep{vanhoeve2022gcdd,KarHoe2026}.

\section{Problem Definition and Arc-Flow Reformulation}
\label{sec:problem}

We next present the vehicle routing problem, the compact mixed-integer programming formulation, and the arc flow reformulations that are derived from dynamic program relaxations.

\subsection{Vehicle Routing with Resource-Constrained Routes}
\label{sec:vrp}

We consider vehicle routing problems of the following form.
Let $C=\{1,\ldots,n\}$ be a set of customers, and let $0$ and $n+1$ denote
copies of the depot corresponding to the start and end of a route,
respectively. We define $V:=C\cup\{0,n+1\}$ as the set of locations, and
\[
E := \{(i,j) \in V \times V : i \neq j, \; i\neq n+1,\; j\neq 0 \}.
\]
as the edge set representing the direct traversals between two locations.
Each customer $i\in C$ has demand $d_i$ and time window $[e_i,l_i]$; we define
$d_0 = d_{n+1} = 0$, $e_0 = e_{n+1} = 0$, and $l_0 = l_{n+1} = T^{\rm max}$ where
$T^{\rm max}$ represents the latest time at which a vehicle can return to the depot.
An unlimited number of identical vehicles is available, each with capacity $Q$.
For each $(i,j)\in E$, let $c_{ij}$ denote the travel cost and $t_{ij}$ the
travel time from $i$ to $j$, where $t_{ij}$ includes any service time at $i$.

A {\em route} is an ordered sequence of locations that begins at $0$ and ends
at $n+1$. Route feasibility is characterized by a set $R$ of resources.
For each resource $r \in R$, we define $q_i^r$ as the resource level at location $i \in V$,
with a given initial value $\bar{q}_0^r$ at the depot. When an edge $(i,j) \in E$ is traversed,
resource $r$ is changed by the (fixed) value $\Delta^r_{ij}$ representing the `consumption' of the resource by the edge. Lastly, each location $i \in V$ has an associated (fixed) interval
$[\ell^r_i,u^r_i]$ representing the admissible range of resource $r \in R$ at that location.
A route is {\em resource feasible} if the implied level for each resource along the route lies within
the admissible intervals of the locations.
The problem is to find a collection of resource-feasible routes that serves every customer exactly once and minimizes the total travel cost.

In this work we consider monotone resources, i.e., having only non-negative consumption, although our framework can also handle non-monotone resources. Specifically, for the capacitated vehicle routing problem with time windows, the capacity and travel time are both monotone resources. For the capacity resource, the initial level is $q_0^{\rm cap} = 0$, the consumption is $\Delta^{\rm cap}_{ij} = d_i$ for $(i,j) \in E$, and the admissible interval is $[\ell_i, u_i] = [0, Q]$ for each $i \in V$.
For the time resource, the initial level is $q_0^{\rm time} = 0$, the consumption is $\Delta^{\rm time}_{ij} = t_{ij}$ for $(i,j) \in E$, and the admissible interval is $[e_i, l_i]$ for each $i \in V$.

\subsection{Compact MILP Formulation}
\label{sec:milp}

We next formulate the vehicle routing problem from Section~\ref{sec:vrp} as a mixed-integer linear program. For each edge $(i,j)\in E$, we introduce a binary variable $x_{ij}$ to indicate whether a
vehicle travels directly from location~$i$ to~$j$. As in Section~\ref{sec:vrp} we use $q^r_i$ as the variable to represent the level of resource $r \in R$ at location $i \in V$.
We define the following compact formulation:
\begin{equation}\label{eq:vrptw}
\min_{x,q} \Big\{ \displaystyle\sum_{(i,j)\in E}c_{ij}x_{ij} : (x, q) \in P ~ \Big\}
\end{equation}
where
\smallskip
\begin{center}
$\begin{array}{rll}
P = \Big\{ ~ (x,q) : & \displaystyle\sum_{i:(i,j) \in E} x_{ij} = 1 & \forall j \in C,\\
    & \displaystyle\sum_{j:(i,j) \in E} x_{ij} = 1 & \forall i \in C,\\
    & q^r_j \geq q^r_i + \Delta^r_{ij} - M^r_{ij}(1-x_{ij})
    & \forall r \in R, \forall (i,j) \in E,\\
    & \ell^r_i \leq q^r_i \leq u^r_i & \forall r \in R, \forall i\in V,\\
    & q^r_0 = \bar{q}^r_0 & \forall r \in R, \\
    & x_{ij} \in \{0,1\} & \forall (i,j)\in E \qquad \Big\}.
\end{array}
$
\end{center}
\medskip

Here $M^r_{ij}$ is a valid big-$M$ constant for resource $r \in R$ and edge $(i,j) \in E$. The first two sets of constraints ensure that every customer has exactly one predecessor and one successor. The remaining constraints propagate the resource levels at each location and enforce resource feasibility along the selected~edges. The resource-propagation constraints also eliminate customer-only subtours, provided that the total increase of at least one resource is strictly positive along every directed cycle~\citep{KULKARNI198558}.

While formulation~\eqref{eq:vrptw} is compact, its continuous relaxation may be weak, in particular because resource feasibility is enforced through
big-$M$ constraints. One approach to strengthening the formulation is to
add valid inequalities in the original variable space, as we will discuss in Section~\ref{sec:vrptw-cuts}. Alternatively, strong relaxations can be obtained by embedding extended formulations that capture the discrete structure from the resource constraints more directly. We next present the latter approach using arc-flow formulations that are derived from dynamic programming models.

\subsection{Dynamic Programming Models for Resource-Constrained Routes}
\label{ss:dpmodels}

We use dynamic programming models to represent sets of routes that reflect feasibility with respect to one or more resources. Existing methods for solving vehicle routing problems use a single dynamic programming model that contains information about all of the resources in the problem, but we consider dynamic programming models that may use a strict subset of the resources, thereby defining relaxations of the set of routes. This permits the use of an ensemble of dynamic programming models that may jointly cover all of the resources.

Let $X$ denote the set of feasible sequences to the vehicle routing problem defined in Section~\ref{sec:problem}. We introduce a family of discrete relaxations of $X$, indexed by $\mathcal H$.
Each relaxation is based on a subset of the resources $R^h \subseteq R$ and, optionally, elementarity of the routes. It is described by a dynamic programming model $\DP^h$.
Each dynamic programming model is defined by states and transitions. A route is encoded in a dynamic program as a sequence of transitions that begins at an initial state and ends at a terminal state. Each state keeps track of the consumption of resources in $R^h$ and optionally enforces elementarity of a route by maintaining a set of visited locations. Each transition corresponds to traversing an edge in $E$ and moves from a starting state to a new state with its resources updated according to the problem description.

We define $\DP^h = (S^h,\tau^h)$ as the dynamic programming model for a subset of resources $R^h \subseteq R$ where $S^h$ is the set of states and
$\tau^h : (S^h \times V) \rightarrow S^h$ is a transition function defined only for admissible state-location pairs.
Each state is defined by a tuple of information $(\{q^{r}\}_{r \in R^h}, i, \NGset)$ where $q^{r}$ is the consumption level of resource $r \in R^h$, $i$ is the current location, and optionally, $\NGset$ is the set of visited locations.
We distinguish an initial state $s^h \in S^h$ corresponding to the start depot and the initial resource levels, and a terminal state $t^h\in S^h$, corresponding to the end depot.
The transition function takes a state and location as input and outputs an updated state based on the consumption of the resources in the problem description for traversing the corresponding edge:
\smallskip
\begin{center}
$\tau^h( (\{q^{r}\}_{r \in R^h}, i, \NGset), j ) = (\{\max\{\ell_j^r,q^{r} + \Delta^{r}_{ij}\}\}_{r \in R^h}, j, \NGset \cup \{j\})$
\end{center}
\smallskip
for locations $j \in V \setminus \NGset$ such that $\max\{\ell_j^r,q^{r} + \Delta^{r}_{ij}\} \leq u^r_j$ for all $r \in R^h$.
We define $X^h$ as the set of routes encoded as sequences of transitions from the initial state to the terminal state through $\DP^h$.

We say that a state $s \in S^h$ {\em dominates} a state $s' \in S^h$ if both represent the same current location and every feasible continuation from $s'$ is also a feasible continuation from $s$ with respect to $\DP^h$. We will use this dominance relation in the next section when constructing the state-transition network. The specific dominance relations used in our arc-flow models are
described in Section~\ref{sec:design}.

In practice, the dynamic programming models defined above may themselves
be too large to represent explicitly. We will therefore typically use relaxations of these models that encode a superset of their routes using fewer states and
transitions. Examples include resource discretization and $ng$-route
relaxations. When such a relaxation is used, it serves as the dynamic programming model to which the arc-flow construction and subsequent aggregation and refinement procedures are applied. Our approach does not depend on the particular dynamic
program used: any such model defines states, transitions, and dominance
relations from which an arc-flow formulation can be constructed. Section~\ref{sec:design} describes the specific relaxations used in our computational implementation.

\subsection{Strengthening the Compact Formulation with Arc-Flow Formulations}
\label{ss:strengthen}

We next strengthen the continuous relaxation of the MILP formulation~\eqref{eq:vrptw} by adding for each discrete relaxation the arc-flow formulation based on the state-transition graph of its dynamic program.
For each relaxation $h \in \mathcal{H}$ with $\DP^h = (S^h, \tau^h)$,
we define the directed network $D^h = (N^h, A^h)$, where $N^h := S^h$.
The arc set is defined as
\medskip

\begin{center}
$
\begin{array}{ll}
A^h := & \big\{ (s, s') \in S^h \times S^h : s' = \tau^h(s, j) \text{~for some~} (s,j) \in \operatorname{dom}(\tau^h) \big\} \cup \\
 & \big\{ (s, s') \in S^h \times S^h :s \text{~dominates~}s' \big\}.
\end{array}
$
\end{center}
\medskip
The first subset defines an arc for each transition, corresponding to a decision in the original problem. The second subset defines arcs based on state dominance; these do not correspond to decisions in the original problem.

For notational convenience we introduce an arc label mapping.
For an arc $a = (s,s') \in A^h$, let
$s = ( \{q^r\}_{r\in R^h}, i, \NGset)$ and $s'=( \{q'^r\}_{r\in R^h}, j, \NGset')$.
We define the arc label mapping by
\[
\ell(a) :=
\left\{
  \begin{array}{ll}
    (i,j), & i \neq j,\\
    0 & i=j.
  \end{array} \right.
\]

The arc-flow formulation of relaxation $h \in \mathcal{H}$ is defined over variables $y^h_a$ representing the flow through arc $a \in A^h$, given by the following constraint set:
\medskip
\begin{center}
\renewcommand{\arraystretch}{1.4}
\begin{tabular}{@{}r@{\;}r@{\qquad}l@{\qquad}r@{}}
$P^h = \bigl\{(x,y^h):$
&
$\displaystyle
\sum_{a\in\delta^+(s)}y_a^h
=
\sum_{a\in\delta^-(s)}y_a^h$
&
$\forall s\in N^h\setminus\{s^h,t^h\},$
&
\refstepcounter{equation}\label{eq:flow}(\theequation)
\\
&
$\displaystyle
x_{ij}=\sum_{\substack{a\in A^h:\\\ell(a)=(i,j)}}y_a^h$
&
$\forall(i,j)\in E,$
&
\refstepcounter{equation}\label{eq:mapping}(\theequation)
\\
&
$y_a^h\geq 0$
&
$\forall a\in A^h\ \bigr\}$
&
\refstepcounter{equation}(\theequation)
\end{tabular}
\end{center}
\medskip
Constraints \eqref{eq:flow} enforce flow conservation at each internal node, and constraints \eqref{eq:mapping} project the flow variables onto the variables in the compact formulation.

Note that we can equivalently formulate the flow conservation constraints \eqref{eq:flow} compactly in matrix form: We let $F^h$ denote the node-arc incidence matrix of $D^h$ with the rows corresponding to $s^h$ and $t^h$ removed, and write $F^h y^h = 0$. We will use this form extensively in later sections.

Observe that these constraints convexify the discrete relaxations $h \in \mathcal{H}$ and can therefore strengthen formulation~\eqref{eq:vrptw}. Specifically, in the special case that only a single exact dynamic program is defined over all resources $R$, this formulation is equivalent to the Dantzig-Wolfe decomposition bound~\citep{sadykov2013column,yamin2026}.
When instead strict subsets of $R$ are used to define the dynamic programs, the formulation is theoretically weaker than the decomposition bound.

We can now define the extended formulation of the original problem, which adds the arc-flow constraints for each discrete relaxation to the original compact formulation:
\begin{equation}\label{eq:ext}
\min \Big\{ \displaystyle\sum_{(i,j) \in E} c_{ij}x_{ij} : (x, q, \{y^h\}_{h \in \mathcal{H}}) \in \Pext ~ \Big\}
\end{equation}
where
\begin{center}
$
\Pext = \big\{ (x, q, \{y^h\}_{h \in \mathcal{H}}) :
(x, q) \in P , \,
(x, y^h) \in P^h ~\forall h \in \mathcal{H}~ \big\}.
$
\end{center}
\medskip
We denote by $z^{\mathrm{ext}}$ the optimal value of the continuous relaxation of the extended formulation~\eqref{eq:ext}. This will be the target bound that we wish to obtain with our method by iterative aggregation and refinement of the network representations.
Recall that in practice we typically do not use the dynamic programs defined in Section~\ref{ss:dpmodels} directly, but instead start from relaxed dynamic programs. In that case, the corresponding arc-flow formulations recover the bound associated with these relaxations, which provides an approximation of the target bound $z^{\mathrm{ext}}$.

\section{Arc-Flow Aggregation}
\label{sec:aggregation}

The arc-flow formulations defined in Section~\ref{sec:problem} may be prohibitively large, motivating the use of smaller aggregated networks. Aggregation reduces the size of a network by merging states, but may weaken the resulting LP relaxation in general. Our approach uses dual information to identify state aggregations that preserve the current LP bound, discussed in this section. We develop a complementary technique of node refinement to strengthen the representation in Section~\ref{sec:refinement}.

Aggregation and refinement are performed separately for each arc-flow
representation. We therefore fix $h\in\mathcal H$ and suppress the
superscript $h$ throughout this section. Thus, $D=(N,A)$ denotes the
corresponding directed graph, $y$ its arc-flow variables, and $F$ its associated node-arc incidence matrix used in the flow conservation constraint definition.

\subsection{Node Aggregation}
\label{sec:node-aggregation}

We obtain a smaller arc-flow formulation by aggregating nodes of
$D=(N,A)$, by partitioning the node set into subsets that will be aggregated.
Let $N'$ be a partition of $N$ such that the source $s \in N$ and terminal $t \in N$ remain singleton nodes. Let $\phi:N\to N'$ map each node to the unique aggregate node containing it. Thus, $\phi(v)$ is the aggregate node containing~$v$.

The partition $N'$ induces an aggregated network $D'=(N',A')$.
Each original arc $a = (v, w) \in A$ induces an arc with tail
$\phi(v)$ and head $\phi(w)$ and retains its label $\ell(a)$.
Arcs with the same aggregate endpoints and label are combined and represented as a single arc, while arcs with different labels are retained as parallel arcs.
We denote the resulting arc mapping by $\psi : A \rightarrow A'$. Thus, $\psi(a)$ is the arc in $A'$ corresponding to $a \in A$. For each $a'\in A'$, we define the corresponding aggregate flow variable
by
\begin{equation}
    y'_{a'}
    :=
    \sum_{\substack{a\in A:\\\psi(a)=a'}} y_a.
\label{eq:aggregate-flow}
\end{equation}
It follows directly that aggregation preserves the projection of the flow variables onto the original variables in equation \eqref{eq:mapping}:
\[
\displaystyle x_{ij}
  = \sum_{\substack{a \in A:\\ \ell(a) = (i,j)}} y_a
  = \sum_{\substack{a' \in A':\\ \ell(a') = (i,j)}} y'_{a'} \qquad \forall (i,j) \in E.
\]

We let $F'$ denote the node-arc incidence matrix associated with $D'$. Observe that each flow-balance equation in $F'y' = 0$ is obtained by summing the flow-balance equations of the original nodes contained in the corresponding aggregate node. This observation will be central to the dual argument in the next subsection.

Applying the above construction to each network $D^h = (N^h,A^h)$, $h \in \mathcal{H}$,
yields an aggregated network $D^{h'} = (N^{h'},A^{h'})$,
where $h'$ denotes the current aggregation associated with $h$, with corresponding
flow variables $y^{h'}$. We denote by $P^{h'}$ the arc-flow constraint set
associated with $D^{h'}$, defined analogously to $P^h$ in
Section~\ref{ss:strengthen}.
We define the aggregated formulation over all networks as
\begin{equation}\label{eq:agg}
\min \Big\{ \displaystyle\sum_{(i,j) \in E} c_{ij}x_{ij} : (x, q, \{ y^{h'}\}_{h \in \mathcal{H}}) \in P^{\mathrm{agg}} ~ \Big\}
\end{equation}
where
\begin{center}
$
P^{\mathrm{agg}} = \big\{ (x, q, \{y^{h'}\}_{h \in \mathcal{H}}) :
(x, q) \in P , \,
(x, y^{h'}) \in P^{h'} ~\forall h \in \mathcal{H}~ \big\}.
$
\end{center}
\medskip
We denote by $z^{\mathrm{agg}}$ the optimal value of the continuous relaxation of the aggregated formulation~\eqref{eq:agg}.

The aggregated formulation is a relaxation of the extended formulation, and therefore $z^{\mathrm{agg}} \leq z^{\mathrm{ext}}$. The inequality may be strict because aggregating nodes may introduce additional feasible flows. We next identify aggregations that preserve the LP bound, using the dual variables associated with the flow-conservation constraints.

\subsection{Dual-Guided Aggregation}
\label{sec:dual-aggregation}
We next identify node aggregations that preserve the LP bound of the current aggregated formulation. Consider an optimal primal-dual solution to the continuous relaxation of $P^{\mathrm{agg}}$, with optimal value $z^{\mathrm{agg}}$.
For the network $D = (N,A)$, let $\pi_v$ denote the dual multiplier associated with the flow-balance equation of node $v \in N \setminus \{s,t\}$. We refer to $\pi_v$ as the potential of node~$v$.

Recall that the flow-balance equation of an aggregate node $v'\in N'$ is obtained by summing the flow-balance equations of the nodes $v$ contained in $v'$. If these equations have the same potential, their contribution to the dual constraints can therefore be represented by assigning this common potential to the aggregate flow-balance equation. This gives the following result.

\begin{theorem}
\label{thm:dual-aggregation}
Let $\pi$ be an optimal dual solution to the continuous relaxation of
$P^{\mathrm{agg}}$, with optimal value $z^{\mathrm{agg}}$. For a network
$D = (N,A)$ in the aggregated formulation, let $N'$ be any partition of $N$
such that nodes are assigned to the same aggregate node only if they have the
same potential. Then replacing $D$ by the aggregated network
$D' = (N',A')$ preserves the optimal LP value $z^{\mathrm{agg}}$.
\end{theorem}

\begin{proof}{Proof}
Replacing $D$ by $D'$ relaxes the primal formulation, so its optimal
value is at most $z^{\mathrm{agg}}$. By construction, the flow-balance equation of
each aggregate node is the sum of the corresponding original
flow-balance equations. Since $\pi$ is constant within each aggregate
node, assigning the common potential to the aggregate flow-balance
equation preserves the corresponding dual constraints and dual
objective value. Hence the aggregated formulation has a dual-feasible
solution of value $z^{\mathrm{agg}}$, so its optimal value is at least
$z^{\mathrm{agg}}$. Thus, the two optimal values are equal.
\Halmos
\end{proof}

Dominance arcs also help explain why equal potentials arise. Such an arc has no routing label, so zero reduced cost implies equal potentials at its endpoints. Chains of tight dominance arcs therefore provide groups of nodes that can be aggregated. In the shortest-path interpretation of the dual, this occurs when relaxing the resource or NG-memory state does not improve the least-cost completion, with arc weights given by the coupling multipliers. This explains a source of aggregation opportunities, without requiring equality in every optimal dual solution; Appendix~\ref{app:lm-equal-potentials} gives the argument in more detail.

Theorem~\ref{thm:dual-aggregation} shows how to further aggregate the networks in $P^{\mathrm{agg}}$ without weakening its LP bound. The initial aggregations, however, may be substantially coarser and yield a weaker LP bound. We therefore disaggregate nodes as needed to recover the target bound $z^{\mathrm{ext}}$. We next develop a refinement procedure for this purpose.

\section{Arc-Flow Refinement}
\label{sec:refinement}

For computational efficiency reasons, we typically initialize our method with networks that are substantially coarser than the original networks defined in Section~\ref{sec:problem}. These aggregations are not required to satisfy the bound-preserving condition of Theorem~\ref{thm:dual-aggregation}, and may therefore yield a weaker LP bound.
The original networks remain available throughout the procedure, however. Starting from an optimal solution of the current aggregated formulation, we construct a restricted LP over the original networks and use its node potentials to identify nodes to disaggregate.

\subsection{Restricted Arc-Flow Formulation}
\label{sec:restricted-af}

Let $(x^*,\{y^{h',*}\}_{h\in\mathcal H})$ be an optimal solution of the
current aggregated LP formulation, with objective value $z^{\mathrm{agg}}$, and let an
optimal dual solution be fixed. For each $h\in\mathcal H$, let
$D^h=(N^h,A^h)$ be the original network and
$D^{h'}=(N^{h'},A^{h'})$ its current aggregation. We denote the corresponding
node and arc mappings by $\phi^h:N^h\rightarrow N^{h'}$ and
$\psi^h:A^h\rightarrow A^{h'}$, respectively.

We construct a restricted set of arcs for each original network with two properties. First, the current optimal solution must remain feasible after applying the same restriction to the aggregated networks. For this purpose, we retain all original arcs that map to aggregate arcs with zero reduced cost
and whose labels correspond to routing edges with positive flow in the current solution, together with dominance arcs.
Second, the restricted original networks must remain feasible. For this purpose, we assume that a feasible integer solution $\bar{x} \in X$ is available, with corresponding integral arc flows $\bar{y}^h$ in each~$D^h$.
Formally, let $\bar r_{a'}^h$ denote the reduced cost of the aggregate flow variable $y_{a'}^{h'}$, and define the restricted arc set
\begin{equation}
\begin{array}{rl}
\widehat A^h := \Bigl\{
& a\in A^h :
\bar y_a^h>0, \, \text{or}\\
&
\bar r_{\psi^h(a)}^{h'}=0
\text{ and }
\bigl(\ell(a)=0 \text{ or } x^*_{\ell(a)}>0\bigr)
\Bigr\}.
\end{array}
\label{eq:restricted-arc-set}
\end{equation}
In particular, $\widehat A^h$ contains all original arcs that map to aggregate arcs carrying positive flow in the current solution.

Let $\widehat{F}^h$ denote the restriction of $F^h$ to the columns indexed by $\widehat{A}^h$. We define the restricted arc-flow constraint set as
\medskip
\begin{center}
$
\begin{array}{rll}
\widehat{P}^h := \Big\{ ~ (x, y^h) : & \widehat{F}^h y^h = 0, \\
& \displaystyle x_{ij} = \sum_{\substack{a \in \widehat{A}^h:\\ \ell(a) = (i,j)}} y^h_a & \qquad \forall (i,j) \in E, \\
& y_a^h \geq 0 & \qquad \forall a \in \widehat{A}^h ~~\Big\}.
\end{array}
$
\end{center}
\medskip

We then define the restricted extended constraint set as
\medskip
\begin{center}
$
\widehat P^{\mathrm{ext}}
:=
\left\{
(x,q,\{y^h\}_{h\in\mathcal H}) :
(x,q)\in P,\;
(x,y^h)\in\widehat P^h
\quad \forall h\in\mathcal H
\right\},
$
\end{center}
\medskip
and let $\widehat z$ denote the optimal value of its continuous relaxation.
This relaxation is feasible because the flows $\bar y^h$ use only arcs in
$\widehat A^h$ for every $h\in\mathcal H$.

For comparison, consider applying the current node aggregation to these same restricted networks. Equivalently, retain in each current aggregated network only the arcs in $\{ \psi^h(a) : a \in \widehat{A}^h \}$.
We refer to the resulting formulation as the restricted aggregated formulation.

\begin{proposition}
\label{prop:restricted-bounds}
The continuous relaxation of the restricted aggregated formulation has optimal value $z^{\mathrm{agg}}$. Moreover,
\begin{equation}
    z^{\mathrm{agg}} \leq z^{\mathrm{ext}} \leq \widehat z.
\label{eq:restricted-bound-relation}
\end{equation}
\end{proposition}

\begin{proof}{Proof}
Restricting the current aggregated formulation can only increase its optimal
value. However, every aggregate arc with positive value
$y_{a'}^{h',*} > 0$ has zero reduced cost. Hence, by~\eqref{eq:restricted-arc-set},
every original arc that maps to such an aggregate arc belongs to
$\widehat{A}^h$. The current optimal solution therefore remains feasible in the
restricted aggregated formulation, which consequently has optimal value
$z^{\mathrm{agg}}$.
The current aggregated formulation is a relaxation of the continuous
relaxation of $P^{\mathrm{ext}}$, so
$z^{\mathrm{agg}} \leq z^{\mathrm{ext}}$. Conversely, the
restricted LP is obtained from the continuous relaxation of $P^{\mathrm{ext}}$
by restricting the arc variables to $\widehat{A}^h$. Hence its feasible region
is contained in that of the continuous relaxation of $P^{\mathrm{ext}}$, which
gives $z^{\mathrm{ext}} \leq \widehat{z}$.
\Halmos
\end{proof}

The restricted LP therefore provides a certificate for whether the current aggregation attains $z^{\mathrm{ext}}$ and, as we show next, the dual information needed to refine it when it does not.

\subsection{Dual-Guided Refinement}
\label{ss:dual-refine}

Let $\widehat{\pi}_v^h$ denote the dual multiplier associated with the
flow-conservation constraint of node $v \in N^h \setminus \{s^h,t^h\}$ in an optimal dual solution to the continuous relaxation over $\widehat{P}^{\mathrm{ext}}$.
Because these multipliers are associated with original-node balance equations, nodes within the same aggregate node may receive different values.

\begin{theorem}
\label{thm:dual-refinement}
If $\widehat{z} > z^{\mathrm{agg}}$, then there exist a network $h \in \mathcal{H}$ and two
original nodes $v, w \in N^h$ such that $\phi^h(v)=\phi^h(w)$ and
$\widehat{\pi}_v^h \neq \widehat{\pi}_w^h$.
\end{theorem}
\begin{proof}{Proof}
Suppose instead that $\widehat{\pi}^h$ is constant on every aggregate node
for every $h \in \mathcal H$. By construction, the flow-balance equation of
each aggregate node is the sum of the corresponding original flow-balance
equations. Assigning their common potential to the aggregate flow-balance
equation therefore preserves the corresponding dual constraints and dual
objective value. Keeping all other dual multipliers unchanged gives a
dual-feasible solution of the restricted aggregated formulation with
objective value $\widehat{z}$. By Proposition~\ref{prop:restricted-bounds}, its continuous relaxation has optimal value $z^{\mathrm{agg}}$. Weak duality then gives
$\widehat{z} \leq z^{\mathrm{agg}}$, contradicting $\widehat{z} > z^{\mathrm{agg}}$. \Halmos
\end{proof}

We use the restricted-LP potentials to refine the current aggregation.
For each aggregate node $v' \in N^{h'}$, we partition the nodes
$v \in N^h$ with $\phi^h(v) = v'$ according to their potentials
$\widehat{\pi}_v^h$. Thus, two original nodes remain in the same
aggregate node if and only if they belong to the same aggregate node
before refinement and have the same restricted-LP potential.
We refer to the resulting partition as the {\em dual-guided partition}.
It is the coarsest refinement of the current partition on which
$\widehat{\pi}^h$ is constant.

\begin{proposition}
\label{prop:refinement-recovers-bound}
Let $\widehat{\pi}$ be an optimal dual solution to the restricted LP, with
optimal value $\widehat{z}$. For each $h \in \mathcal H$, let the
dual-guided partition refine the current partition according to
$\widehat{\pi}^h$. Then replacing the original restricted networks by the
networks induced by these partitions preserves the optimal LP value
$\widehat{z}$.
\end{proposition}
\begin{proof}{Proof}
The refined formulation is an aggregation of the continuous relaxation over
$\widehat{P}^{\mathrm{ext}}$, so its LP value is at most $\widehat{z}$.
By construction, $\widehat{\pi}^h$ is constant within each aggregate node.
The same argument as in Theorem~\ref{thm:dual-aggregation} therefore gives a
dual-feasible solution of the refined formulation with objective value
$\widehat{z}$. Weak duality gives the reverse inequality, and the two values
are equal. \Halmos
\end{proof}

The above results provide a characterization of the refinement step. If $\widehat{z} = z^{\mathrm{agg}}$, the current aggregation already attains the restricted LP bound and no refinement is required. If $\widehat{z} > z^{\mathrm{agg}}$, Theorem~\ref{thm:dual-refinement} guarantees that the dual-guided partition refines at least one aggregate node, while Proposition~\ref{prop:refinement-recovers-bound} shows that the networks induced by these partitions attain the LP value $\widehat{z}$ on the restricted arc sets.

\section{Reduced Cost-Based Arc Removal}
\label{sec:arc-removal}

While the node aggregation can substantially reduce the size of the arc flow networks, we can additionally remove arcs that cannot belong to an improving integer solution based on a reduced cost-based argument that is common for arc-flow formulations~\citep{delima_et_al_2023,karahalios2023column}. We extend this approach to a lower bound that is based on using multiple arc flow formulations for the same compact formulation.

Suppose that an integer feasible solution $\bar x\in X$ with objective value $\bar z=c^\top\bar x$ is available. Consider an optimal primal-dual solution of the LP over $P^{\mathrm{agg}}$, with objective value $z^{\mathrm{agg}}$.
As before, for each $h\in\mathcal H$, let $D^h=(N^h,A^h)$ denote the original network and $D^{h'}$ its current aggregation. We assign to each original node $v\in N^h$ the potential of its aggregate node $\phi^h(v)$ and compute the reduced costs $\bar{r}_a^h$ of the original arcs using these potentials and the dual multipliers of the linking constraints.

For an arc $a = (v,w) \in A^{h}$, let
$d_{s^{h}v}^{h}$ denote the minimum reduced cost of a path from
$s^{h}$ to $v$, and let $d_{w t^{h}}^{h}$ denote the minimum reduced
cost of a path from $w$ to $t^{h}$. We define
\begin{equation}
    \rho_{a}^{h}
    :=
    d_{s^{h}v}^{h}
    + \bar{r}_{a}^{h}
    + d_{wt^{h}}^{h}.
\end{equation}
as a lower bound on the reduced cost of an $s^{h}$-$t^{h}$ path containing $a$.
The coupling through the common variables $x$ allows this reduced-cost bound to be strengthened across the other networks. For each routing edge $(i,j) \in E$, we define
\begin{equation}
    \rho_{ij}^{h}
    :=
    \min\left\{
        \rho_{a}^{h} :
        a \in A^{h},\ \ell(a) = (i,j)
    \right\}
\label{eq:decision-path-reduced-cost}
\end{equation}
as the minimum reduced cost of an $s^{h}$-$t^{h}$ path in $D^{h}$ that uses routing edge $(i,j)$.
Lastly, for each $a \in A^{h}$, we define
\begin{equation}
    L_{a}^{h}
    :=
    \begin{cases}
    z^{\mathrm{agg}} + \rho_{a}^{h}
    + \displaystyle\sum_{\widetilde{h} \in \mathcal{H} \setminus \{h\}}
        \rho_{ij}^{\widetilde{h}},
        & \text{if } \ell(a) = (i,j), \\[2ex]
    z^{\mathrm{agg}} + \rho_{a}^{h},
        & \text{if } \ell(a) = 0.
    \end{cases}
\label{eq:arc-removal-bound}
\end{equation}
By construction, $L_{a}^{h}$ is a lower bound on the objective value of any integer feasible solution using~$a$. This gives the following test for removing arcs.
\begin{proposition} \label{prop:reduced-cost}
For $h \in \mathcal{H}$ and $a \in A^{h}$, if
$L_{a}^{h} \geq \bar{z}$, then no integer feasible solution with
objective value strictly less than $\bar{z}$ can use $a$. Hence, $a$ can be removed from $A^{h}$ without eliminating an improving solution.
\end{proposition}
The forward and backward distances used to compute $\rho_{a}^{h}$ can be obtained by shortest-path calculations on the original networks. Since the reduced costs are nonnegative by dual feasibility, standard shortest-path algorithms can be used.
Arcs satisfying the condition of Proposition~\ref{prop:reduced-cost} are removed from the original networks, and the aggregated networks are updated accordingly through $\psi^h$, i.e., an aggregate arc is removed when none of its corresponding original arcs remains. Lastly, if all arcs labeled $(i,j)$ are removed from any one network, the routing edge $(i,j)$ can also be removed from the compact formulation and all other networks. The removal test can be repeated whenever the aggregated LP solution or incumbent value changes.

\section{Integrating Problem-Specific Cutting Planes}
\label{sec:vrptw-cuts}

An advantage of retaining the compact formulation \eqref{eq:vrptw} is that
cutting planes can be added directly in the original $x$-space, without
modifying the arc-flow representations. Therefore, to further strengthen the continuous relaxation of~\eqref{eq:vrptw}, we use two families of inequalities:
{\em minimum-route inequalities} and {\em local-area partition cuts}. Both exploit feasibility information that is not captured by the individual route
relaxations. Minimum-route inequalities impose a lower bound on the number of routes required to serve a customer set, while local-area partition cuts further restrict how these customers can be served through feasible route segments.

The first family generalizes standard capacity and $k$-path inequalities for
vehicle routing. For a customer set $S\subseteq C$, such inequalities bound
the number of routes that must cross the boundary of $S$. The classical
capacity argument gives the lower bound
\[
    \left\lceil
        \frac{\sum_{i\in S} d_i}{Q}
    \right\rceil
\]
on the number of vehicles required to serve $S$. We strengthen this argument
by accounting for both capacity and time-window feasibility.

For a set of customers $S\subseteq C$, let $k(S)$ denote the minimum number
of feasible routes required to serve all customers in $S$. In particular,
$k(S)=1$ if the customers in $S$ can be served by a single capacity- and
time-feasible route.
Let $\mathcal F(S)$ contain the nonempty subsets of $S$ that can each be served by a single feasible depot-to-depot route. Then
{
\begin{equation}
 \begin{aligned}
    k(\emptyset)&=0,\\
    k(S)&=\min_{T\in\mathcal F(S)}\{1+k(S\setminus T)\}
    \quad(S\ne\emptyset).
 \end{aligned}
    \label{eq:lm-route-count}
\end{equation}
}
Appendix~\ref{app:lm-cuts} describes the single-route feasibility test and evaluation of this recurrence. For the inequalities below, we assume that deleting customers from a feasible route and shortcutting the resulting gaps preserves resource feasibility. We refer to the resulting inequalities
as minimum-route inequalities.

Consider $\widehat C\subseteq C$ and any $S$ such that
$\widehat C\subseteq S\subseteq C$. Since serving $\widehat C$ requires at
least $k(\widehat C)$ routes, and each such route must leave $S$, the
following minimum-route inequality is valid:
\begin{equation}
    \sum_{(i,j)\in\delta^+(S)} x_{ij}
    \geq k(\widehat C),
\label{eq:minimum-route}
\end{equation}
where
\[
    \delta^+(S)
    :=
    \{(i,j)\in E:i\in S,\;j\notin S\}.
\]
For a fixed $\widehat C$, the most violated inequality in
\eqref{eq:minimum-route} can be found by a minimum-cut computation with
capacities given by the current values of $x$. We restrict attention to
inclusion-minimal sets $\widehat C$ for which removing any customer decreases
$k(\widehat C)$.

The local-area partition cuts strengthen this counting argument by imposing additional consistency on the routing decisions within a small customer set. They are obtained through a Benders separation problem. Given a fractional solution $x$, the subproblem determines whether the values of the routing variables associated with the customer set can be represented by feasible route segments. When this is not possible, an optimal dual solution of the subproblem yields a Benders cut in the original $x$-variables. Since the separation problem may admit multiple optimal dual solutions, we apply the Pareto-optimal cut selection procedure of \cite{magnanti1981accelerating} to select a stronger representative from this
family.

Consider a small customer set $C' \subseteq C$ and let $\kappa := k(C')$.
If an integer solution serves $C'$ using exactly $\kappa$ route segments, then
the customers in $C'$ can be partitioned into $\kappa$ nonempty sets, each
traversed by a feasible contiguous route segment. The possible segments for
each set are obtained from a dynamic program that enumerates the
nondominated feasible orderings of its customers.

Among equal-cost segments, the dominance rule uses the lexicographic order of their customer visit lists. This ordering is preserved under a common extension at either end, as required by the optimality-preservation argument. The labels, dominance rule, and recursion defining this efficient frontier are given in Appendix~\ref{app:lm-frontier}. By restricting attention to nondominated segments, the resulting inequalities may exclude feasible integer solutions, but under the stated dominance and tie-breaking rules, preserve a common optimal solution across the entire cut family (Proposition~\ref{prop:lm-la-preservation}).

Since at least $\kappa$ route segments are required to serve $C'$, an integer solution with fewer than $\kappa+1$ boundary crossings must use exactly $\kappa$ segments. Accordingly, given a fractional solution $x$ for which the total flow leaving $C'$ is less than $\kappa+1$,
the LA-$\kappa$ separation LP assigns nonnegative weights to partitions into $\kappa$ blocks and to frontier segments serving those blocks. Internal arc flows and boundary entry and exit flows limit these weights. A slack variable measures the portion of the required partition weight $\kappa+1-\sum_{(i,j)\in\delta^+(C')}x_{ij}$ that cannot be supported. A positive optimal slack yields a separating inequality through the LP dual. Appendix~\ref{app:lm-separation} gives the separation model and cut derivation. In our implementation, we use LA-$1$ and LA-$2$ cuts.

Both cut families therefore strengthen the compact formulation directly,
while the arc-flow networks continue to impose the discrete route
relaxations defined in Section~\ref{sec:vrptw-relaxations}. The cuts and the
arc-flow representations can consequently be generated and maintained
independently while interacting through the common variables $x$.

\section{Algorithmic Framework}
\label{sec:algorithm}

We next integrate the methods developed in the previous sections into an algorithmic framework, establish finite convergence of the aggregation and refinement procedure, and describe the construction of the final MILP.

\subsection{Dual-Guided Aggregation and Refinement}

We combine the dual-guided aggregation and refinement procedures in an iterative procedure presented as Algorithm~\ref{alg:aggregation-refinement} below. We assume that for each of the original networks $D^h$ for $h \in \mathcal{H}$ an initial aggregation $D^{h'}$ is given. We also assume that an initial feasible solution $\bar{x}$ with associated integer flow variables $\bar{y}^h$ is given.

\begin{algorithm}[t]
\caption{Dual-guided aggregation and refinement}
\label{alg:aggregation-refinement}
\begin{algorithmic}[1]
\setlength{\baselineskip}{1.4em}
\REQUIRE Original networks $D^h$ and initial aggregations $D^{h'}$ for each $h \in \mathcal{H}$, and initial integer feasible solution $(\bar{x},\{\bar{y}^h\}_{h\in\mathcal H})$
\STATE Set $z^{\mathrm{last}}\leftarrow -\infty$
\REPEAT
    \STATE Solve the aggregated LP, obtaining $(x^*,\{y^{h',*}\}_{h \in \mathcal H})$, value $z^{\mathrm{agg}}$, and node potentials $\{\pi^{h'}\}_{h \in \mathcal H}$
    \IF{$z^{\mathrm{agg}} > z^{\mathrm{last}}$}
        \STATE Apply dual-guided aggregation using $\{\pi^{h'}\}_{h \in \mathcal H}$
        \STATE Set $z^{\mathrm{last}}\leftarrow z^{\mathrm{agg}}$
    \ENDIF
    \STATE Construct restricted arc set $\widehat A^h \subseteq A^h$ for each $h\in\mathcal H$ %
    \STATE Solve the restricted LP, obtaining optimal value $\widehat{z}$ and node
    potentials $\{\widehat\pi^{h}\}_{h \in \mathcal H}$
    \IF{$\widehat{z} > z^{\mathrm{agg}}$}
        \STATE Apply dual-guided refinement using $\{\widehat\pi^{h}\}_{h \in \mathcal H}$
    \ENDIF
\UNTIL{$\widehat{z} = z^{\mathrm{agg}}$}
\end{algorithmic}
\end{algorithm}

Algorithm~\ref{alg:aggregation-refinement} proceeds by solving the LP relaxation with respect to the current aggregated networks. If the optimal solution improves the last bound $z^{\rm last}$, we apply the dual-guided aggregation procedure from Section~\ref{sec:dual-aggregation} and update $z^{\rm last}$. Observe that $z^{\rm last}$ is used to avoid repeated aggregation at the same LP value, thereby preventing cycling.
We then construct the restricted arc set for each network and solve the associated LP over the restricted original networks. If the resulting solution improves the aggregated LP bound, we apply the dual-guided refinement procedure from Section~\ref{ss:dual-refine}. We repeat the process until the aggregated and restricted LP bounds coincide.

\begin{theorem}
\label{thm:finite-convergence}
Algorithm~\ref{alg:aggregation-refinement} terminates after finitely many
iterations with LP value $z^{\mathrm{ext}}$.
\end{theorem}

\begin{proof}{Proof}
At every iteration, Proposition~\ref{prop:restricted-bounds} gives
$
z^{\mathrm{agg}} \leq z^{\mathrm{ext}} \leq \widehat{z}.
$
Hence, if $\widehat{z} = z^{\mathrm{agg}}$, then $z^{\mathrm{agg}} = z^{\mathrm{ext}}$.
It remains to show finite termination. Refinement can only increase the
current LP value, while dual-guided aggregation preserves it by
Theorem~\ref{thm:dual-aggregation}. Since the original networks are finite, there are only finitely many node partitions and therefore only finitely many possible LP values. The condition $z^{\mathrm{agg}} > z^{\mathrm{last}}$ thus allows dual-guided aggregation only finitely many times.

After the final aggregation step, every iteration with $\widehat{z} > z^{\mathrm{agg}}$
strictly refines at least one aggregate node by
Theorem~\ref{thm:dual-refinement}. Only finitely many such refinements are
possible. At full refinement, the current formulation is the original
arc-flow formulation and therefore has value $z^{\mathrm{agg}}$. Its optimal
arc flows are contained in the restricted arc sets by construction, so they
are feasible for the restricted LP. Hence
$\widehat{z} \leq z^{\mathrm{ext}}$, and
Proposition~\ref{prop:restricted-bounds} gives
$\widehat{z} = z^{\mathrm{ext}} = z^{\mathrm{agg}}$. The algorithm therefore terminates. \Halmos
\end{proof}

\subsection{Integrated Algorithm and Final MILP Construction}

Theorem~\ref{thm:finite-convergence} establishes finite convergence to
$z^{\mathrm{ext}}$, the target bound defined by the original arc-flow
formulations. By separating the cutting planes from
Section~\ref{sec:vrptw-cuts}, we can additionally capture interactions
or feasibility conditions that are not represented by the individual
arc-flow formulations and can therefore strengthen the relaxation
beyond $z^{\mathrm{ext}}$. We integrate these cutting planes, together
with the reduced cost-based arc removal from
Section~\ref{sec:arc-removal}, into
Algorithm~\ref{alg:aggregation-refinement}.

Algorithm~\ref{alg:aggregation-refinement-cuts} summarizes the
modifications to the original procedure. Lines~7a and 7b introduce
reduced cost-based arc removal and cutting plane separation before
constructing the restricted arc sets. The remaining steps are
unchanged, while the modified termination condition in line~11'
requires that no violated cutting planes are found and
$\widehat z=z^{\mathrm{agg}}$.

\begin{algorithm}[t]
\caption{Extensions to Algorithm~\ref{alg:aggregation-refinement}}
\label{alg:aggregation-refinement-cuts}
\begin{algorithmic}
\setlength{\baselineskip}{1.4em}
\STATE \textbf{Insert before line 7:}
\item[7a:] Apply reduced cost-based arc removal
\item[7b:] Separate violated cutting planes and add them to the formulation
\STATE \textbf{Replace line 11 by:}
\item[11':] \textbf{until} no violated cutting planes are found
and $\widehat z=z^{\mathrm{agg}}$
\end{algorithmic}
\end{algorithm}

After termination, we apply a final dual-guided aggregation to
further reduce the formulation size without weakening its LP bound.
We then simplify the aggregated networks by removing redundant arcs
whose effect can be reproduced by an alternative path involving
dominance arcs and the same routing decision. This operation
preserves the feasible projection onto the original variables.

The removal criterion and the pre-compression of the restricted LP are detailed in Appendix~\ref{sec:restricted-construction}.
Finally, we apply a final round of reduced cost-based arc removal as described in Section~\ref{sec:arc-removal}, including arcs retained during
refinement to preserve feasibility of the restricted LP. The
resulting aggregated networks, together with the compact formulation
and generated cutting planes, define the MILP passed to the solver.

\section{Arc-Flow Model Design}
\label{sec:design}
\label{sec:vrptw} %

We next describe our design choices in defining the discrete relaxations for the vehicle routing problem, the dynamic programs, and the initial aggregated networks. We also discuss how we apply reduced cost-based arc filtering to further reduce the size of the networks.

\subsection{Discrete Relaxations and Dynamic Programs}
\label{sec:vrptw-relaxations}

The vehicle routing problem consists of two types of resources, vehicle capacity and travel time (to enforce the time windows), as well as elementarity of the routes. In principle, one can define a separate relaxation based on every subset of these three elements. In our implementation we found that two relaxations performed best: one that is based on vehicle capacity, and one that combines travel time and elementarity.
In fact, for the latter we apply the ng-route relaxation~\citep{baldacci2011new} and time-discretization to obtain a relaxation that is more tractable in practice.
Both relaxations are defined through a dynamic programming transition model, described next.

\smallskip

\noindent\emph{Capacity.~}
The capacity relaxation, which we will refer to as CAP, retains only the vehicle-capacity restriction. A
state is given by $(i,q)$, where $i\in V$ is the current location and $q$ is
the remaining vehicle capacity. The initial state is $(0,Q)$. From a state
$(i,q)$, the next location may be a customer $j\in C$ with $d_j\leq q$, in
which case the transition is
\[
(i,q) \xrightarrow{j} (j,q-d_j).
\]
Alternatively, the route may return to the end depot $n+1$.
The capacity state space contains every reachable integer capacity level. We note that in the benchmarks we considered all capacities and demands are integer valued.

\smallskip

\noindent\emph{Elementarity + Time.~}
The second relaxation combines time-window feasibility with ng-route elementarity, which we will refer to as NGT.
For each location $i\in C$, let $C_i^{\mathrm{ng}}\subseteq C$ denote its NG neighborhood~\citep{baldacci2011new}.
In addition, for each location $i \in V$, we partition its time window $[e_i,l_i]$
into contiguous intervals of width $\delta$, with the first and last
intervals truncated as necessary.
A state is given by $(i,I,M)$, where $i \in V$ is the current location, $I$ is one of the time intervals associated with $i$, and $M \subseteq C_i^{\mathrm{ng}}$ is the NG memory.
A transition to customer $j$ is permitted if $j \notin M$ and $\max\{e_j,\underline I+t_{ij}\}\leq l_j$, where $\underline I$ denotes the lower endpoint of $I$.
Let $I'$ denote the time interval at $j$ containing $\max\{e_j,\underline I+t_{ij}\}$. After visiting $j$, the memory is updated as $M' = C_j^{\mathrm{ng}} \cap (M \cup \{i\})$.
The resulting transition is
\[
    (i,I,M) \xrightarrow{j} (j,I',M').
\]
In the implementation, the arrival time defining $I'$ is also tightened by a bound computed from the tracked route history; Appendix~\ref{app:lm-graph-details} specifies this bound and its computation.

\noindent\emph{Dominance.~} The dominance relations for these dynamic programs follow directly from
their state representations. For the capacity relaxation, a state
$(i,q)$ dominates $(i,q')$ if $q \geq q'$.
For the ngt relaxation, $(i,I,M)$ dominates $(i,I',M')$ if
$\underline{I} \leq \underline{I}'$ and $M \subseteq M'$. In constructing the
corresponding networks, we exploit transitivity of these relations and
include dominance arcs only between immediately dominating states.
The sparse arc construction, including the resource bound used for NGT transitions, is specified in Appendix~\ref{app:lm-graph-details}.

For each $h \in \{\mathrm{CAP}, \mathrm{NGT}\}$, the states
and transitions from their dynamic program define the network $D^h=(N^h,A^h)$ and its
associated arc-flow formulation $P^h$. These dynamic programs are
tractable relaxations of the models defined in Section~\ref{ss:dpmodels},
and together they induce a relaxation of $P^{\mathrm{ext}}$. The tractability is particularly
relevant for the refinement step, which relies on these original networks
to define the restricted arc-flow formulation.

\subsection{Initial Network Aggregations}
\label{sec:vrptw-initial-aggregation}

Recall that our algorithm solves linear programs over aggregated
and restricted networks. Because the state-transition networks defined above
are typically too large to be optimized over in reasonable time, we initialize
our method with aggregated networks that provide a coarse approximation.

For the capacity relaxation, all states $(i,q)$ associated with the same
location $i$ are initially aggregated into a single node, so the initial
network has one node per location.
For the NGT relaxation, we construct the initial aggregation by partitioning the time intervals at each location into consecutive groups (in our implementation three intervals) and disregarding the NG memory set.
Thus, states $(i,I,M)$ whose time intervals belong to the same
group are represented by a single aggregate node, irrespective of their NG
memory $M$. The refinement procedure may subsequently disaggregate these nodes according to the underlying time intervals and NG-memory states.

\section{Computational Experiments}
\label{sec:experiments}

The computational study examines whether dual-guided aggregation and
refinement makes strong arc-flow formulations practical. We first evaluate
the quality of the root bounds and the effort required to obtain them,
including the contribution of cuts. We then examine the size of the
formulations constructed by the procedure. Finally, we test their use for
integer optimization through branch and bound. We refer to the implementation
as Dual-Guided Aggregation and Refinement (DGAR).

\subsection{Experimental Setup}
\label{sec:exp-setup}

We use the 168 Solomon VRPTW instances with 25, 50, and 100 customers
\citep{solomon1987algorithms}. Each size contains 17 clustered (C), 23 random
(R), and 16 mixed (RC) instances. The objective is total travel cost with an
unrestricted fleet and no fixed vehicle cost; distances are truncated to one
decimal place. We retain the size--family breakdown throughout the analysis
and distinguish RC1 from RC2 where this helps explain the mixed-family results.

DGAR combines aggregation, refinement, reduced-cost removal, and separation in the original
routing variables. It uses the NGT relaxation described in
Section~\ref{sec:vrptw-relaxations}, a fixed NG neighborhood parameter of five,
and the additional capacity relaxation on C instances. Aggregation and
refinement adapt the state-space representation while the NG neighborhoods
remain fixed. Time discretization has
spacing 10 for C, 1 for R1 and RC1, and 5 for R2 and RC2. After root
strengthening, the final explicit model is passed to Gurobi with its routing
variables declared binary. The graphs are fixed during the subsequent search.

Appendix~\ref{app:lm-frontier-audit} documents the frontier rule used in these runs and the separate comparison with the theoretical rule.

We compare DGAR with the compact formulation and three route-based methods.
The \emph{compact}
formulation supplies the baseline for the improvement in relaxation strength
and for the final integer solve. \emph{BP-NG5} and \emph{BP-NG9} provide
route-based bounds from branch-and-price implementations with
fixed NG neighborhoods of five and nine. \emph{RouteOpt}
\citep{you2025routeopt} is a branch-price-and-cut solver that provides a further
route-based benchmark. We measure its root relaxation at the last completed
exact-pricing checkpoint before the first cut-separation round, including
NG-memory augmentation and the standard fleet lower-bound row. The comparison
assesses the bound quality delivered by these configurations, with DGAR's cuts
enabled and each method's strengthening mechanisms as described above.
For integer optimization, all five methods are evaluated through their
complete search procedures; RouteOpt continues with branch-price-and-cut.

All experiments were performed on the same AMD Ryzen 9 5900XT workstation
with 16 physical cores, 32 logical processors, and 64~GB of RAM, running
Windows 11. The DGAR and compact runs use Gurobi 13.0.0.
Instances were solved sequentially. DGAR and the
compact formulation used up to 32 solver threads; RouteOpt used 32 LP threads,
with serial pricing and separation. DGAR was given 1,800 seconds for root
strengthening, followed by a 3,600-second limit for the final integer solve.
The compact formulation and RouteOpt had a 3,600-second solve limit.
A separate paired single-thread experiment measures the contribution of
DGAR's cuts (Appendix~\ref{app:lm-cut-ablation}).

For every method we compute the root gap using the same benchmark reference
objective $U_i$ and the method's lower bound $L_{im}$:
\begin{equation}
 g_{im}=100\,\frac{U_i-L_{im}}{U_i}.
 \label{eq:exp-common-gap}
\end{equation}
The reference values are the reported optima for this travel-cost benchmark.
We report arithmetic means of per-instance gaps. Two lower bounds are treated
as tied when their absolute difference is at most 0.05 objective units,
half the reporting increment of RouteOpt's one-decimal bounds.
The gap profiles count the instances for
which each method obtains a valid completed root bound within a given gap
threshold. They cover all 168 instances; a run without a valid completed
bound contributes no step. Table~\ref{tab:exp-main} reports each method's
mean gap over its valid completed roots, with coverage counts for BP-NG5 and
BP-NG9. Coverage details are given in Appendix~\ref{app:exp-audit}.

DGAR's reduced-cost removal receives the reference upper bound $U_i+10^{-4}$
and a separate feasible backup solution. Obtaining the reference bound is outside the measured effort. The compact
runs do not receive a supplied incumbent, so the integer comparison includes
the benefit of DGAR's supplied upper-bound information.

\subsection{Bound Quality and Root Effort}
\label{sec:exp-bound-quality}\label{sec:exp-bounds}

DGAR substantially strengthens the compact relaxation. Across all 168 instances, it reduces
the mean root gap from 23.03\% to 0.90\%, a reduction of 96.10\%. Its gap is
at most 1.00\% on 122 instances and at most 5.00\% on 164. The improvement appears
at every size: DGAR's mean gaps are 0.46\%, 0.82\%, and 1.41\% at 25, 50,
and 100 customers, whereas the compact means remain between 20.00\% and 25.00\%.
The improvement is therefore sustained as the number of customers grows,
while the formulation remains explicit.

\begin{table}[t]
\centering
\caption{Root-bound quality by size and family. Mean gaps use each method's
valid completed roots; parentheses give their number for BP-NG5 and BP-NG9.
The other three methods complete all $n$ roots in each row.
G/T/R counts stronger DGAR bounds, ties,
and stronger RouteOpt bounds.}
\label{tab:exp-main}
\small
\setlength{\tabcolsep}{3.5pt}
\begin{tabular}{rlrrrrrrr}
\toprule
& & & \multicolumn{5}{c}{Mean gap (\%)} & \\
\cmidrule(lr){4-8}
Size & Family & $n$ & Compact & BP-NG5 & BP-NG9 & DGAR & RouteOpt & G/T/R \\
\midrule
25 & C & 17 & 13.73 & 0.03 (17) & 0.03 (17) & 0.01 & 0.00 & 0/16/1 \\
25 & R & 23 & 19.57 & 0.78 (23) & 0.65 (23) & 0.31 & 0.58 & 12/9/2 \\
25 & RC & 16 & 28.69 & 2.07 (14) & 0.79 (15) & 1.15 & 0.74 & 1/12/3 \\
\addlinespace[2pt]
50 & C & 17 & 11.77 & 0.00 (16) & 0.00 (15) & 0.00 & 0.00 & 0/16/1 \\
50 & R & 23 & 23.54 & 1.10 (16) & 0.97 (17) & 1.15 & 1.01 & 12/2/9 \\
50 & RC & 16 & 39.02 & 6.57 (12) & 5.51 (12) & 1.23 & 4.19 & 7/5/4 \\
\addlinespace[2pt]
100 & C & 17 & 11.82 & 0.00 (11) & 0.00 (10) & 0.06 & 0.00 & 0/14/3 \\
100 & R & 23 & 25.03 & 1.15 (13) & 0.90 (12) & 1.96 & 1.13 & 4/1/18 \\
100 & RC & 16 & 36.52 & 2.87 (9) & 2.66 (9) & 2.06 & 1.90 & 9/0/7 \\
\midrule
All & & 168 & 23.03 & 1.41 (131) & 1.11 (130) & 0.90 & 1.02 & 45/75/48 \\

\bottomrule
\end{tabular}
\end{table}

The route-based comparisons put this gain in perspective. On all 168
instances, RouteOpt has a mean gap of 1.02\%, versus 0.90\%
for DGAR, and reaches the 1.00\% threshold on 117 instances versus 122. The
pairwise counts are 45 stronger DGAR bounds, 75 within tolerance, and 48
stronger RouteOpt bounds. The size and family breakdown in
Table~\ref{tab:exp-main} explains the differences behind these aggregate results.

The spatial structure of the instances matters. Both methods have nearly
closed bounds on C. On R with 25 customers, DGAR's
mean gap is 0.31\%, compared with 0.58\%; at 100 customers the corresponding
means are 1.96\% and 1.13\%. The clearest advantage occurs on RC1: over its
24 instances, DGAR has 16 stronger bounds and eight within tolerance, with
mean gaps of 0.42\% and 4.22\%. At 50 customers, its RC1 mean gap is just
0.06\%, compared with 8.27\%. RC2 exhibits the opposite ordering, with
means of 2.54\% and 0.33\%.

\begin{figure}[t]
\centering
\includegraphics[width=\linewidth]{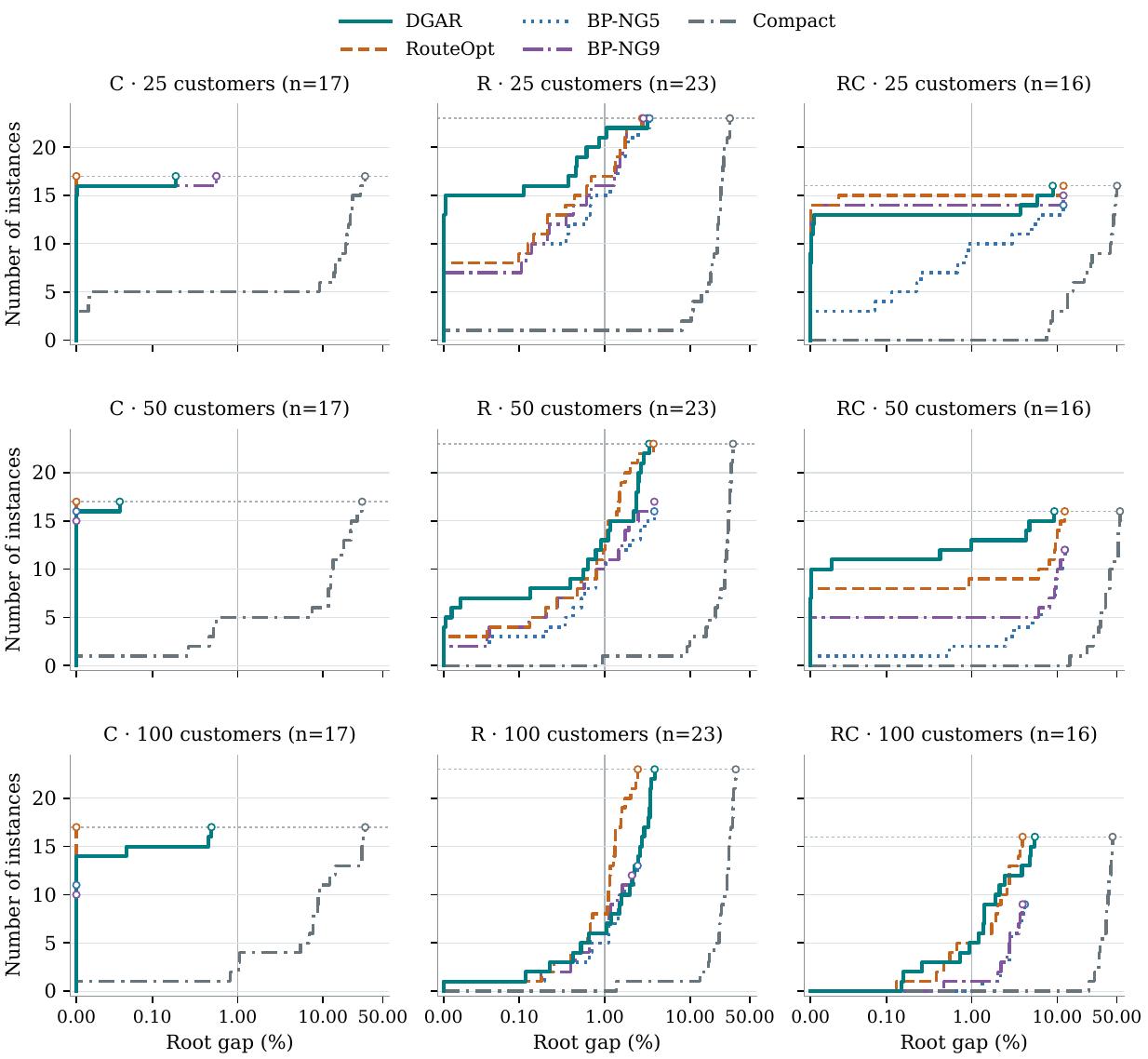}
\caption{Cumulative numbers of instances within each root-gap threshold,
for all five methods over the full 168-instance benchmark. Each panel gives
its instance count $n$, also marked by a horizontal gray line. Every valid
completed bound contributes one step. Each curve ends at its largest observed
gap, with a circle marking the endpoint and the number of completed bounds.
The horizontal axis is linear up to 0.10\% and logarithmic thereafter;
the vertical line marks a gap of 1.00\%.}
\label{fig:exp-nine}
\end{figure}

Figure~\ref{fig:exp-nine} compares all five methods over the full benchmark.
At a gap threshold of 1.00\%, BP-NG5 and BP-NG9 reach 86 and 94 instances,
respectively, compared with 122 for DGAR. Their curves also reveal the
difference in root completion: BP-NG5 and BP-NG9 provide valid completed
bounds on 131 and 130 instances, while DGAR, compact, and RouteOpt complete
all 168. The profiles therefore show both the strength of the bounds and
the number of instances on which they are obtained.

On the 128 instances with valid completed bounds for all five methods,
DGAR's mean gap is 0.34\%, compared with 1.42\% for BP-NG5 and 1.12\%
for BP-NG9. This common set is easier for DGAR than the full benchmark,
on which its mean gap is 0.90\%. The paired one-thread ablation in
Appendix~\ref{app:lm-cut-ablation} reduces the mean gap from 1.86\% to
0.83\% when cuts are enabled, with the largest improvement on RC. These
measurements concern the combined cut families within the procedure.

\paragraph{Root computational effort.}
\label{sec:exp-time}

A strong relaxation is useful only if its construction cost is controlled.
In the 32-thread runs, DGAR completes the root procedure on all 168
instances. Median wall times, including setup and graph construction, are
3.64, 29.90, and 167.01 seconds at the three sizes. The corresponding
90th percentiles are 37.59, 175.78, and 596.47 seconds. These times quantify
the work needed to turn the weak compact relaxation into the strengthened
explicit model.

\begin{table}[t]
\centering
\caption{Root effort at 32 threads, in seconds, on all 56 instances per size.
DGAR includes program setup and root strengthening; RouteOpt uses elapsed time
through completed exact pricing; Compact times its separate plain LP
solve. The last column is DGAR's 90th percentile.}
\label{tab:exp-time}
\small
\begin{tabular}{rrrrr}
\toprule
& \multicolumn{3}{c}{Median root time (s)} & DGAR 90th \\
\cmidrule(lr){2-4}
Size & Compact & DGAR & RouteOpt & percentile (s) \\
\midrule
25 & 0.02 & 3.64 & 0.42 & 37.59 \\
50 & 0.09 & 29.90 & 1.39 & 175.78 \\
100 & 0.48 & 167.01 & 7.33 & 596.47 \\

\bottomrule
\end{tabular}
\end{table}

Table~\ref{tab:exp-time} also shows the cost of the alternative bound
computations. The compact LP is inexpensive but leaves substantially larger
gaps. DGAR generally requires more root work than RouteOpt. This additional effort builds the explicit formulation used for integer
optimization. Appendix~\ref{app:exp-scatter} gives the paired gap and time
comparisons for all 168 instances.

\begin{figure}[t]
\centering
\includegraphics[width=\linewidth]{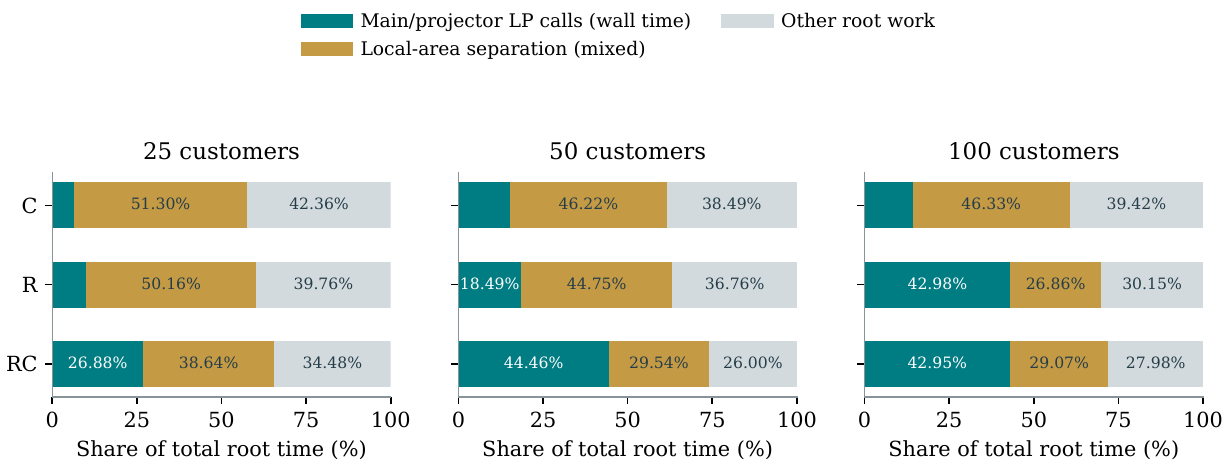}
\caption{DGAR root-time decomposition at 32 threads, using sums within each
size--family group. Main/projector LP calls include solver-interface wall
time. Local-area separation contains both solver and non-solver work. Other
root work is the residual, including setup and graph/model operations. The
three categories sum to the measured root time.}
\label{fig:exp-time-split}
\end{figure}

To understand DGAR's cost, Figure~\ref{fig:exp-time-split} separates three
non-overlapping portions of its measured root time. Calls solving the main
lower-bound and projector LPs account for 38.85\% of the total. Local-area
separation accounts for 31.02\%, and the remaining root work for 30.13\%.
The latter includes graph operations, model construction, and setup.
Gurobi's own solution time for the main and projector LPs accounts for
38.51\% of total root time, close to the wall time of these calls. Local-area
separation also solves auxiliary LPs, so its cost combines solver time with
the construction and processing of the separation problems. The decomposition
identifies work on all three stages as relevant to reducing root effort.

\subsection{Formulation Size}
\label{sec:exp-formulation-size}\label{sec:exp-aggregation}

Strong arc-flow bounds become practical when they can be represented by
models of manageable size. We therefore follow the formulation through four
stages of the 32-thread runs: the original model constructed from the
unaggregated networks, the first aggregated lower-bound LP, the final
bound-producing LP, and the model passed to the integer solver after closing
cleanup. Table~\ref{tab:exp-stages} reports total variable counts at each
stage, including routing, arc-flow, and auxiliary variables, and identifies
the arc-flow component of the final model. All counts
precede Gurobi presolve. The original model is constructed but not optimized
in these runs; it provides the size reference for the subsequent reductions.

\begin{table}[t]
\centering
\caption{Evolution of formulation size at 32 threads. Entries are medians
over the $n$ instances in each size--family group; ``All'' pools the three
families within a size. The four stage columns count all variables.
Final $y$ counts the arc-flow variables within the final MILP.
Final counts are measured after closing cleanup and before Gurobi presolve.}
\label{tab:exp-stages}
\footnotesize
\setlength{\tabcolsep}{3pt}
\begin{tabular}{rlrrrrrr}
\toprule
Size & Family & $n$ & Original model & First LP & Final LP & Final MILP & Final $y$ \\
\midrule
25 & C & 17 & 183,639.00 & 5,217.00 & 1,907.00 & 250.00 & 151.00 \\
25 & R & 23 & 434,798.00 & 16,721.00 & 1,342.00 & 233.00 & 119.00 \\
25 & RC & 16 & 331,369.00 & 11,523.50 & 1,994.00 & 215.50 & 105.00 \\
25 & All & 56 & 337,170.00 & 11,261.00 & 1,804.00 & 233.00 & 131.00 \\
\addlinespace[3pt]
50 & C & 17 & 397,019.00 & 17,720.00 & 2,730.00 & 438.00 & 260.00 \\
50 & R & 23 & 895,925.00 & 56,894.00 & 4,452.00 & 1,273.00 & 893.00 \\
50 & RC & 16 & 696,991.50 & 42,236.00 & 3,609.50 & 485.00 & 253.50 \\
50 & All & 56 & 738,360.50 & 40,739.50 & 3,983.50 & 487.50 & 278.50 \\
\addlinespace[3pt]
100 & C & 17 & 999,623.00 & 64,211.00 & 5,933.00 & 851.00 & 501.00 \\
100 & R & 23 & 1,978,478.00 & 223,322.00 & 47,696.00 & 12,635.00 & 9,070.00 \\
100 & RC & 16 & 1,491,123.50 & 158,937.50 & 35,578.50 & 9,114.50 & 6,692.50 \\
100 & All & 56 & 1,612,628.00 & 153,137.50 & 21,316.00 & 5,619.00 & 3,991.50 \\

\bottomrule
\end{tabular}
\end{table}

The initial aggregation makes the first LP substantially smaller than the
original model at every size. For 100 customers, the median falls from
1,612,628.00 to 153,137.50 variables. Subsequent iterations strengthen the
relaxation while reducing the median to 21,316.00 variables at the final
lower-bound LP. Over these same 56 instances, the mean gap decreases from
9.58\% at the first LP to 1.41\% at the final LP. The corresponding mean
gaps decrease from 10.17\% to 0.46\% at 25 customers and from 12.51\%
to 0.82\% at 50 customers. Thus, the stronger bounds are obtained with
smaller working LPs. The final LP has fewer variables than the first on
154 of the 168 instances; refinement can also enlarge the representation
to strengthen the relaxation.

Closing cleanup further reduces the models passed to integer optimization,
whose median total variable counts are 233.00, 487.50, and 5,619.00 at the
three sizes. These changes reflect the combined effects of aggregation,
refinement, cuts, and pruning. The original construction also includes
artificial depot connectors that are removed before the integer solve,
so the full reduction includes this cleanup. The family breakdown shows
that the final 100-customer models are smallest on C, with a median of
851.00 variables, compared with 12,635.00 on R and 9,114.50 on RC.

\begin{figure}[t]
\centering
\includegraphics[width=\linewidth]{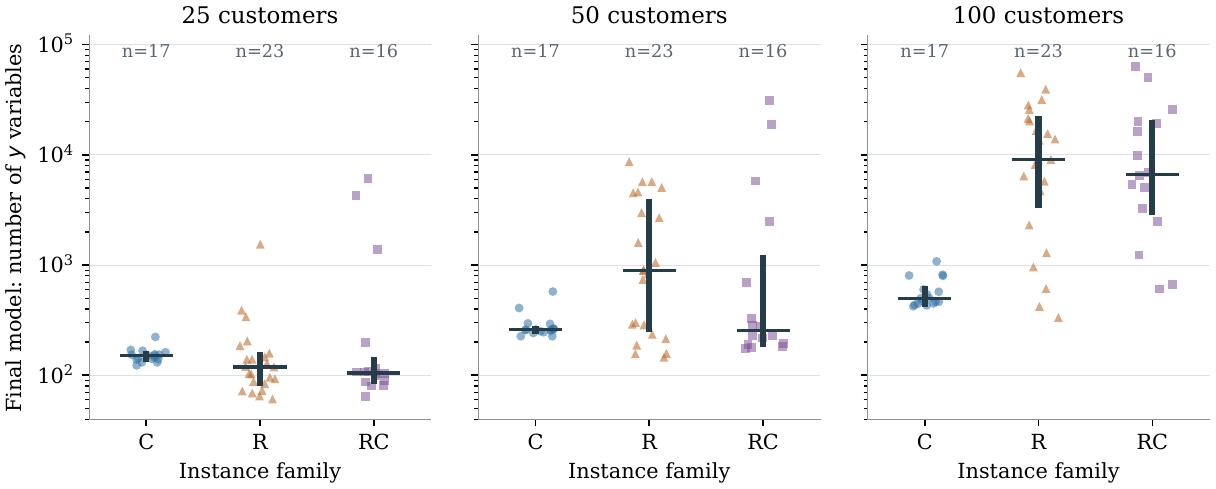}
\caption{Final arc-flow variable counts at 32 threads, for all 168 DGAR
instances. Each point is an instance; horizontal bars mark medians and
vertical bars interquartile ranges. Horizontal displacement separates
observations within categories, and the vertical axis shows actual counts
on a logarithmic scale.}
\label{fig:exp-y-size}
\end{figure}

At 25, 50, and 100 customers, the median final $y$ counts are 131.00, 278.50,
and 3,991.50 (Table~\ref{tab:exp-stages}). Thus, a typical final model
contains hundreds to a few thousand arc-flow variables. At 100 customers
the family medians are 501.00 for C, 9,070.00 for R, and 6,692.50 for RC; the
full range is shown in Figure~\ref{fig:exp-y-size}. These counts describe the
models after DGAR's final pruning and before Gurobi presolve. Together with
the preceding LP counts, they show how the procedure turns the original
networks into strengthened explicit models for integer optimization.
Final row counts are reported in Appendix~\ref{app:exp-time}.
Supplementary measurements of pre-compression in the restricted LP are
reported in Appendix~\ref{sec:restricted-construction}.

\FloatBarrier
\subsection{Integer Optimization}
\label{sec:exp-integer}\label{sec:exp-performance}

As a final test, we pass the strengthened formulation to a general-purpose
MILP solver and compare the ability of all five methods to prove optimality.
The routing variables are declared binary, and the graphs remain fixed
during the search.

Gurobi solves the fixed DGAR model with up to 32 threads, a 3,600-second
MILP limit, \texttt{MIPFocus=3}, and \texttt{Cuts=2}. Its ordinary presolve,
cutting planes, heuristics, and branching operate on the fixed DGAR formulation.
All 168 runs return feasible
integer solutions. They establish optimality at the reference objective on
146 instances: all 56 with 25 customers, 54 with 50, and 36 with 100. Of these,
119 terminate during presolve or after at most one explored node. This count
includes Gurobi's root processing and cuts.

\begin{table}[t]
\centering
\caption{Instances proved optimal within 3,600 seconds of recorded elapsed
 time. Each entry requires a reported solver proof and agreement with the
 reference objective within 0.05 objective units. All $n$ instances remain
 in each comparison. Configurations are described in Section~\ref{sec:exp-setup}.}
\label{tab:exp-integer}
\small
\begin{tabular}{rrrrrrr}
\toprule
Size & $n$ & Compact & BP-NG5 & BP-NG9 & DGAR & RouteOpt \\
\midrule
25 & 56 & 52 & 47 & 51 & 56 & 56 \\
50 & 56 & 25 & 27 & 29 & 54 & 54 \\
100 & 56 & 19 & 13 & 12 & 36 & 49 \\
\midrule
All & 168 & 96 & 87 & 92 & 146 & 159 \\

\bottomrule
\end{tabular}
\end{table}

Table~\ref{tab:exp-integer} compares the five methods within 3,600 seconds.
DGAR proves \mbox{optimality} on 146 instances,
compared with 96 for Compact, 87 for \mbox{BP-NG5}, 92 for \mbox{BP-NG9}, and 159 for
RouteOpt. DGAR and RouteOpt each solve all 56 instances with 25 customers
and 54 with 50 customers; their counts at 100 customers are 36 and 49.
The explicit formulation thus supports substantially more proofs than the
compact and fixed-neighborhood branch-and-price formulations, while
RouteOpt solves more of the largest instances.

Relative to Compact, DGAR gains 50 proofs and retains all 96 of its successes.
The largest increase is at 50 customers, from 25 to 54. Its mean final
incumbent gap over all 168 instances is 0.12\%, compared with 0.65\% for
Compact. These results show that the strengthened formulation can repay
its construction cost through more effective integer search, with the
upper-bound initialization described in Section~\ref{sec:exp-setup}.

\begin{figure}[t]
\centering
\includegraphics[width=\linewidth]{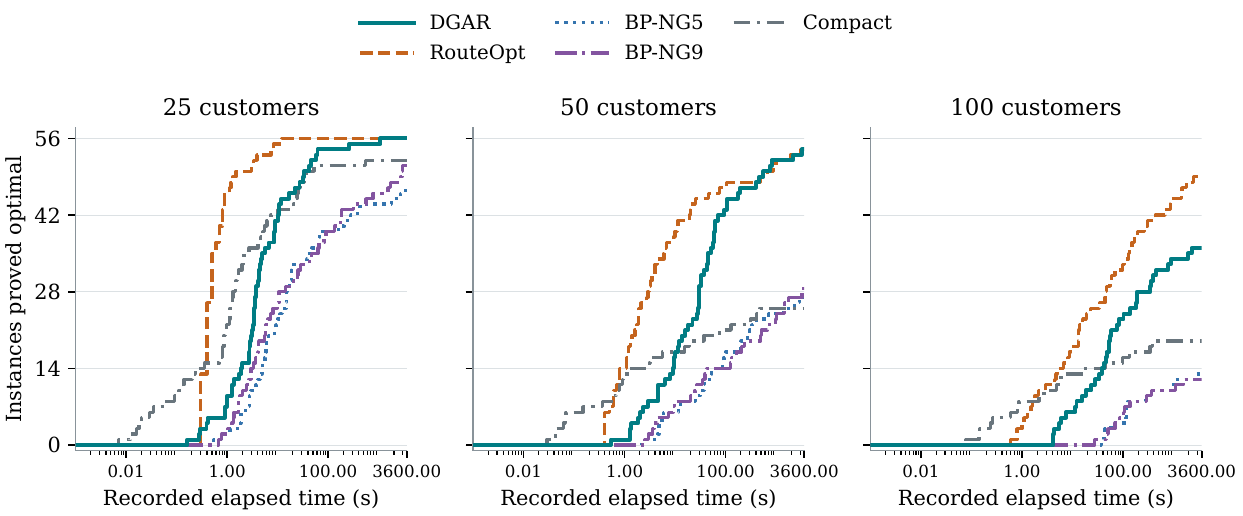}
\caption{Instances proved optimal by recorded elapsed time for all five
methods. Each panel contains 56 instances and uses the same 3,600-second
reporting window. DGAR includes setup, root strengthening, and the final MILP;
BP-NG5, BP-NG9, and RouteOpt include root and tree search; Compact uses MILP
optimization time. Runs without a qualifying proof contribute no step.}
\label{fig:exp-integer-profile}
\end{figure}

Figure~\ref{fig:exp-integer-profile} shows how the five methods accumulate
proofs over time. DGAR's elapsed time includes setup, root strengthening, and
the final MILP; the branch-and-price methods include their root and tree
search, while Compact uses its MILP optimization time. Each curve retains
all 56 instances of its size class and rises only when an optimality proof
is completed. The comparison therefore accounts for the root investment
needed before DGAR can begin integer search.

\FloatBarrier
Branch and bound currently operates on the fixed formulation constructed
by DGAR at the root. It benefits from this strengthened model, while DGAR's
aggregation and refinement mechanisms are not yet used within the search
tree. A tighter integration could adapt the formulation to branching
decisions and strengthen bounds during search. Exploring this integration
is a direction for future work, with the aim of achieving state-of-the-art
performance through compact arc-flow formulations.
\FloatBarrier

\section{Conclusion}
\label{sec:conclusion}

We developed a general methodology for reducing the size of arc-flow
formulations while retaining the strength of their LP relaxations. Starting
from a family of arc-flow representations of discrete relaxations, we showed
that nodes with equal flow-balance dual multipliers can be aggregated without
changing the LP bound. When an aggregation is too coarse, a restricted LP on
the original networks guides node disaggregation to strengthen
the formulation. Combining these operations yielded an iterative procedure
that converges finitely to the LP bound defined by the original
arc-flow representations.

The framework allows strong arc-flow relaxations to be represented by
smaller explicit formulations. By retaining the original routing variables,
problem-specific cutting planes can be incorporated directly, while
reduced-cost information can be used to remove arcs that cannot participate
in an improving integer solution. The resulting formulation combines these
strengthening techniques in a model that can be solved by a general-purpose
MILP solver.

We applied the methodology to the VRPTW using two arc-flow relaxations, one based on vehicle capacity and the other combining time windows with relaxed elementarity. We also developed and integrated two families of cutting planes to further strengthen the formulation. The computational results demonstrate that strong LP bounds can be obtained using substantially aggregated networks, yielding formulations that are orders of magnitude smaller than the original arc-flow formulations. The resulting root bounds were competitive with those obtained by the branch-and-price and branch-cut-and-price methods considered in our experiments. Moreover, these bounds were obtained without column generation or pricing, while retaining an explicit formulation that can be solved directly by a general-purpose MILP solver for integer optimization.

\bibliographystyle{informs2014} %
\bibliography{mip-arcflow} %

\clearpage

\begin{APPENDICES}
\section{Implementation Details}

{
\subsection{Equal Potentials along Dominance Arcs}
\label{app:lm-equal-potentials}
Choose the signs of the flow-balance and coupling multipliers so that an
arc $(v,w)$ has reduced cost $\bar c_{vw}=b_{vw}+\pi_w-\pi_v$, where
$b_{vw}$ is the coupling multiplier for its routing label and is zero for
a dominance arc. Then dual feasibility gives
\[
 \pi_v\le b_{vw}+\pi_w.
\]
A tight dominance arc therefore has $\pi_v=\pi_w$, and this equality
extends along any chain of tight dominance arcs. Complementary slackness
makes every dominance arc carrying positive primal flow tight; tightness
can also occur at zero flow.

For fixed non-flow multipliers, normalize the terminal potential to zero.
On nodes that can reach the terminal, with no negative-weight directed
cycle, let $d_v$ be the minimum weight of a path from $v$ to the terminal.
Every feasible potential satisfies $\pi_v\le d_v$, and $d$ satisfies the
Bellman inequalities with equality on a minimum-weight continuation.
Flow-balance rows have zero right-hand sides, so replacing internal
potentials by this feasible maximal choice leaves the dual objective
unchanged. In networks with no source or terminal flow-balance rows,
these two boundary potentials stay fixed; the Bellman equality is asserted
only for internal nodes. A dominance arc is tight whenever its endpoint
nodes have the same minimum completion weight. Small changes in the
resource or NG-memory state can leave this weight unchanged, explaining
why equal potentials arise. The argument identifies a mechanism, not a
guarantee of large groups or equality in every dual optimum.

\subsection{Arc Construction and Resource Bounds}
\label{app:lm-graph-details}
The abstract construction in Section~\ref{sec:problem} defines the network before
sparse representation and history-dependent tightening. Dominance arcs need
only connect immediate successors in the dominance order: their transitive
closure provides the other dominance transitions. A routing arc can be
omitted when an immediate dominance successor of its tail has an arc with
the same label to the same head. The omitted arc is then represented by a
dominance arc followed by that routing arc. The head must be the same node,
not merely a node at the same location. Source and terminal arcs in the
protected set are retained.

For the NGT network, the destination memory is
$M'=C_j^{\rm ng}\cap(M\cup\{i\})$, as in
Section~\ref{sec:vrptw-relaxations}. The implementation also tightens the
arrival resource using the history retained at $j$. Put
$S=C_j^{\rm ng}\cap M$, with the penultimate location $i$ named separately.
Let $\mathcal Q(S,i,j)$ be the feasible sequences from the depot that serve
exactly $S$ before visiting $i$ and then $j$. In the cumulative-resource
convention of this paper, define
\[
 q^{r,\min}_{S,i,j}
   =\min_{\sigma\in\mathcal Q(S,i,j)}q_j^r(\sigma).
\]
The minima are componentwise; different resources may attain them on
different sequences. For time, the destination interval contains
\[
 \max\{e_j,\underline I+t_{ij},q^{{\rm time},\min}_{S,i,j}\}
\]
rather than just $\max\{e_j,\underline I+t_{ij}\}$.
An infeasible destination or an empty sequence class yields no transition.
The direct depot case is evaluated directly. Under the shortcutting
assumption in Section~\ref{sec:vrptw-cuts}, deleting untracked customers
cannot increase resource consumption, so this minimum is a lower bound on
the cumulative resource of any full route with that retained history.
Thus the tightening preserves those routes.

The quantity denoted $r^+$ in the remaining-resource description is the
same bound with the sign convention reversed: set $\rho^r=-q^r$, so
$r^{r,+}_{S,i,j}=-q^{r,\min}_{S,i,j}$. It is an upper bound on remaining
resource, not a new resource variable. It can be computed by the frontier
recursion of Appendix~\ref{app:lm-frontier}, starting at the fixed depot
resource. For each retained sequence to $i$, append $(i,j)$, test
feasibility, and minimize the resulting arrival resource. Cost and the
visit-order tie-break are unnecessary for this resource-bound computation;
resource dominance suffices. The visit-order rule is needed for the LA
optimality argument instead.

Finally, the implementation may declare routing variables free in its LPs:
nonnegativity follows from the coupling equations, and customer degree
constraints imply the upper bound of one on customer-incident edges.
This statement does not bound a direct depot-to-terminal edge representing
an empty route. Such an edge must be omitted, fixed to zero, or retain its
explicit bound. The compact formulation retains its stated variable bounds.

Declaring the routing variables free also affects the dual information used
for arc removal in Section~\ref{sec:arc-removal}. At an LP optimum, a free
variable $x_{ij}$ has zero reduced cost. Its objective coefficient is
therefore accounted for by the coupling multipliers and the multipliers of
the remaining constraints. With explicit variable bounds, part of this
information may instead be carried by bound multipliers. In particular,
a nonnegative variable at zero can have a positive reduced cost that is
not included in the arc-path quantities $\rho_a^h$ and $\rho_{ij}^h$.
The free-variable representation places this contribution in the coupling
duals used to price the network arcs. It can thus strengthen the
arc-removal certificates without changing the feasible region or the LP
bound. This is the reason for omitting the redundant routing-variable
bounds in the implemented LPs.
}

\subsection{Implementation of the Restricted LP}
\label{sec:restricted-construction}

The restricted arc sets in Section~\ref{sec:restricted-af} are constructed
from the support of a feasible reference flow and the zero-reduced-cost arcs
of the current aggregated LP. The latter can be identified directly from
the current primal-dual solution. The original arcs corresponding to the
retained aggregate arcs are identified through the mapping $\psi^h$. As the
aggregation changes, the restricted sets can therefore be updated without
rebuilding the original networks.

{
\paragraph{Protected arcs.}
The reference flow retained in~\eqref{eq:restricted-arc-set} certifies
feasibility of the restricted LP, but need not place every network node
on a source--terminal path. The implementation imposes this additional
property by retaining a protected set $A^h_{\rm safe}$ of source-incident
and terminal-incident arcs. More precisely, it uses the enlarged set
\[
 \widehat A^h_{\rm impl}
   =\widehat A^h\cup A^h_{\rm safe},
\]
where $\widehat A^h$ is defined in~\eqref{eq:restricted-arc-set}.
The protected routing labels also remain available: fixing $x_{ij}=0$
would prevent every arc with label $(i,j)$ from carrying flow, even if
those arcs were retained.

In the implemented construction, the protected arcs provide source and
terminal connections for every retained node. For a general network,
retaining existing incident arcs alone does not establish this property;
source--terminal reachability must also hold. This connectivity requirement
is stronger than the feasibility requirement in
Section~\ref{sec:restricted-af} and ensures that each retained node has a
route continuation for the potential calculations.

The enlargement leaves the bounds and refinement arguments unchanged.
It retains the feasible reference flow, and its aggregated counterpart
still has value $z^{\mathrm{agg}}$: it contains the current aggregate
optimum and is a restriction of the full aggregated formulation.
If $\widehat z_{\rm impl}$ denotes the enlarged restricted-LP value, then
\[
 z^{\mathrm{agg}}\le z^{\mathrm{ext}}
 \le\widehat z_{\rm impl}\le\widehat z.
\]
Thus Proposition~\ref{prop:restricted-bounds} and the argument of
Theorem~\ref{thm:dual-refinement} apply to this restricted set as well.
Protected arcs are exempt from reduced-cost removal during refinement.
The final MILP may remove them when the usual removal certificate applies.

\paragraph{Pre-compression of the restricted LP.}
This operation precedes the restricted LP solve and is different from
aggregation using an already computed dual solution.
Let $\mathcal Y_v$ be the set of pairs $(\ell(a),w)$ reached from node $v$
by zero or more dominance arcs followed by one routing-labelled arc $a$
ending at $w$. Starting with singleton groups $M$, define
\[
 \operatorname{sig}(v,M)=
 \bigl(\operatorname{location}(v),
       \{(b,M_w):(b,w)\in\mathcal Y_v\}\bigr).
\]
Keep source and terminal nodes separate. Merge nodes with equal signatures,
recompute the signatures, and repeat until no further merge occurs.
Thus nodes in a group have the same labelled continuation options to the
same groups, including the options reached through dominance arcs.
The signature describes outgoing continuation constraints, not identical
rows of the incidence matrix: incoming arcs may differ.

For fixed coupling multipliers, these nodes admit the same maximal feasible
potential, given by the least-weight completion to the terminal as in
Appendix~\ref{app:lm-equal-potentials}. Consequently the restricted LP admits
an optimal dual constant on each merged group. Aggregation relaxes the
primal while this dual preserves its value, so pre-compression leaves the
restricted optimum unchanged. After solving, lift each group potential to
its original nodes and use these potentials in Section~\ref{sec:refinement}.
On an acyclic network, a reverse topological pass suffices; otherwise the
signatures are updated until stable.

\paragraph{Removal of redundant arcs.}
An arc $(v,w)$ can be removed if the retained network contains a path from
$v$ to $w$ with the same routing label contribution: dominance arcs, one
arc with the same nonzero label, and further dominance arcs. For a
dominance arc, the replacement path contains only dominance arcs.
Rerouting any flow on the removed arc along this path preserves flow
balance and every coupling equation, and hence the projection onto $x$.
Identical parallel arcs can likewise be replaced by a single flow variable.
A deletion must retain its replacement path, or a certified replacement of
that path; checking witnesses in a graph and then deleting them all without
this condition would not suffice. Protected arcs are kept during refinement.
This structural removal is separate from reduced-cost removal, which uses
an incumbent objective value.
}

{
\paragraph{Size reduction before solving the restricted LP.}
Table~\ref{tab:restricted-lp-precompression} reports supplementary measurements
of this pre-compression on 168 Solomon instances. The restricted networks
retain all original nodes before aggregation, although only a subset of their
arcs is retained. We compare the node counts before and after aggregation
for two solves in each run: the median-ranked solve and the solve with the
largest post-compression node count. For 100 customers, the median count
falls from 213\,879 nodes before aggregation to 1\,960 for the typical solve
and 2\,317 for the largest solve. The corresponding median per-instance
reduction factors are 116 and 96. Across the three sizes, the median reduction
factor for each run's largest solve ranges from 96 to 133.
These supplementary measurements come from the separate restricted-LP study;
they are not identified as the 32-thread runs of Section~\ref{sec:experiments}.
These counts measure the flow-balance part of the restricted LP, rather than
the total variable counts of the main formulation in
Section~\ref{sec:exp-formulation-size}. They quantify a reduction in LP
dimensions, not a measured runtime speedup.

\begin{table}[!htbp]
\centering
\caption{Restricted-LP node counts before and after pre-compression.
There are 56 instances at each size; All pools the 168 instances.
Within each run, the typical solve is the median-ranked solve by
post-compression node count, and the largest solve has the greatest such
count. All entries are medians over instances. Reduction factors are
computed per instance before taking the median and therefore need not
equal ratios of the displayed median node counts.}
\label{tab:restricted-lp-precompression}
\footnotesize
\setlength{\tabcolsep}{4pt}
\begin{tabular}{rrrrrr}
\toprule
 & Before & \multicolumn{2}{c}{Typical solve} & \multicolumn{2}{c}{Largest solve} \\
\cmidrule(lr){3-4}\cmidrule(lr){5-6}
Size & Nodes & Nodes after & Reduction ($\times$) & Nodes after & Reduction ($\times$) \\
\midrule
25 & 54\,706 & 192 & 219 & 309 & 133 \\
50 & 111\,618 & 712 & 137 & 896 & 109 \\
100 & 213\,879 & 1\,960 & 116 & 2\,317 & 96 \\
\midrule
All & 96\,930 & 736 & 135 & 941 & 112 \\
\bottomrule
\end{tabular}
\end{table}

}
\FloatBarrier

{
\subsection{Restricted-Action Warm Start for the LP}
\label{app:lm-action-warmstart}
This procedure concerns the solution of an LP at a fixed aggregation;
it is separate from the restricted LP used to compute refinement potentials.
Its seed consists of routing edges that carried flow in earlier iterations,
all source- and terminal-incident routing edges, and the five cheapest
outgoing edges of each customer. First solve with all other routing
variables fixed to zero. Then release these temporary restrictions and
solve the full LP, using the first basis as a warm start for primal simplex.
Only the second optimum is used as the lower bound. If nonnegativity of
$x$ is supplied by the arc-flow coupling, an upper bound of zero suffices
to fix a routing variable to zero. Each such fixing also sets every
associated arc-flow variable to zero in every network.

The final unrestricted solve is essential: the restricted optimum alone
is an upper bound on the full LP optimum, not the lower-bound certificate
required by the algorithm. No pricing-based stopping certificate is used
by this procedure. In the reported implementation it was enabled for the
R and RC single-thread runs and disabled for C instances. It was not used
in the 32-thread results of Section~\ref{sec:experiments}; those root times
must not be credited with the separate single-thread warm-start savings.
}

{
\section{Construction and Separation of the Cutting Planes}
\label{app:lm-cuts}

This appendix supplies the dynamic programs and separation model used to define
the two cut families in Section~\ref{sec:vrptw-cuts}. We use the locations,
resources, and routing variables of Section~\ref{sec:problem}.

The frontier below uses a visit-order tie-break to obtain the closure
property needed by the LA-cut proof. Appendix~\ref{app:lm-frontier-audit}
distinguishes this rule from the implementation used in the reported runs
and summarizes the available comparison experiments.

\subsection{Feasible Segments and Efficient Frontiers}
\label{app:lm-frontier}

For the VRPTW, a resource level propagates along $(i,j)\in E$ as
$q_j^r=\max\{\ell_j^r,q_i^r+\Delta_{ij}^r\}$, provided that $q_j^r\le u_j^r$.
The maximum permits waiting for the time resource; for capacity its lower
bound is zero. We describe the frontier using these cumulative resource
levels, consistently with the compact formulation.

Consider a segment $\sigma=(v_0,\ldots,v_m)$ with $m\ge1$. Its summary consists
of a cost $c(\sigma)$ and resource vectors $D(\sigma)$, $E(\sigma)$, and
$L(\sigma)$, defined componentwise by
\begin{align}
 c(\sigma)&=\sum_{h=0}^{m-1}c_{v_hv_{h+1}}, \label{eq:lm-summary-a}\\
 D^r(\sigma)&=\sum_{h=0}^{m-1}\Delta^r_{v_hv_{h+1}}, \nonumber\\
 E^r(\sigma)&=\max_{0\le h\le m}
 \left\{\ell^r_{v_h}+\sum_{a=h}^{m-1}\Delta^r_{v_av_{a+1}}\right\}, \nonumber\\
 L^r(\sigma)&=\min_{0\le h\le m}
 \left\{u^r_{v_h}-\sum_{a=0}^{h-1}\Delta^r_{v_av_{a+1}}\right\}.
 \label{eq:lm-summary-b}
\end{align}
Empty sums are zero. For a feasible segment, $L^r$ is the largest admissible
resource level at its first location, and the resource level at the last
location is $\max\{q_{v_0}^r+D^r,E^r\}$. Feasibility also requires that the
intermediate bounds are respected. The recursion below enforces this
condition when extending a segment.

Suppose that a feasible segment $\sigma$ starts at $w$, and prepend an edge
$(u,w)$. The new summary, denoted by primes, is
\begin{align}
 c'&=c_{uw}+c, \label{eq:lm-prepend-a}\\
 D^{r\prime}&=\Delta^r_{uw}+D^r,\nonumber\\
 E^{r\prime}&=\max\{\ell_u^r+D^{r\prime},E^r\}, \nonumber\\
 L^{r\prime}&=\min\{u_u^r,L^r-\Delta^r_{uw}\}.
 \label{eq:lm-prepend-b}
\end{align}
The extension is feasible if and only if $L^{r\prime}\ge\ell_u^r$ for every
resource $r$. The base case is a single edge $(u,v)$, retained when
$\max\{\ell_v^r,\ell_u^r+\Delta_{uv}^r\}\le u_v^r$ for every $r$.

For distinct endpoints $u,v$ and a customer set $T$ containing neither
endpoint, consider all feasible elementary segments from $u$ to $v$ whose
interior customers are exactly $T$. Write $\omega(\sigma)$ for the ordered
list of visited locations. Within this class, a segment $\gamma$ dominates
$\sigma$ when
\begin{equation}
 \begin{aligned}
 c(\gamma)&\le c(\sigma),& D(\gamma)&\le D(\sigma),\\
 E(\gamma)&\le E(\sigma),& L(\gamma)&\ge L(\sigma),
 \end{aligned}
 \label{eq:lm-dominance}
\end{equation}
and $(c(\gamma),\omega(\gamma))$ is lexicographically smaller than
$(c(\sigma),\omega(\sigma))$. The vector comparisons are componentwise.
The ordering of visit lists breaks cost ties consistently under a common
extension at either end. In particular, a strict improvement in a resource
summary alone is not sufficient for discarding an equal-cost segment.

Let $\mathcal E(u,T,v)$ denote the segments in this class that are not
dominated. These frontiers can be computed in increasing order of $|T|$. Let
$\mathcal A(u,T,v)$ contain all feasible extensions $(u,w)\oplus\sigma$
with $w\in T$ and $\sigma\in\mathcal E(w,T\setminus\{w\},v)$. Then
\begin{align}
 \mathcal E(u,\emptyset,v)
 &=\begin{cases}\{(u,v)\},&\text{if $(u,v)$ is feasible},\\
 \emptyset,&\text{otherwise},\end{cases}\label{eq:lm-frontier-base}\\
 \mathcal E(u,T,v)
 &=\operatorname{ND}\bigl(\mathcal A(u,T,v)\bigr).
 \label{eq:lm-frontier-rec}
\end{align}
Here $\oplus$ denotes prepending and $\operatorname{ND}$ removes dominated
segments using~\eqref{eq:lm-dominance} and its lexicographic condition.

\begin{lemma}
\label{lem:lm-frontier}
Every contiguous subsegment of a frontier segment, containing at least one
edge, belongs to the frontier of its own class.
Recursion~\eqref{eq:lm-frontier-rec} returns the complete frontier.
\end{lemma}
\begin{proof}{Proof}
If a subsegment were dominated, replacing it by a dominating segment would
preserve feasibility and weakly improve the summary of the complete segment.
The common prefix and suffix preserve the strict lexicographic ordering of
cost and visit list. The replacement would therefore dominate the original
segment, a contradiction. In particular, deleting the first edge of a
frontier member leaves a frontier member. Induction on $|T|$ now shows that
every frontier member occurs in $\mathcal A(u,T,v)$.
Every other member of this set has a dominating frontier member, since
the class is finite. Removing dominated members therefore returns exactly
the frontier. \Halmos
\end{proof}

\subsection{Single-Route Feasibility and Minimum-Route Separation}
\label{app:lm-minimum}

For a depot-to-depot route, fix the initial resource levels by setting
$\ell_0^r=u_0^r=\bar q_0^r$ during the frontier computation. A customer set
$T$ is single-route feasible if and only if
\begin{equation}
 \mathcal E(0,T,n+1)\ne\emptyset.
 \label{eq:lm-single-route}
\end{equation}
Computing these tests and~\eqref{eq:lm-route-count} in increasing subset size
determines $k(S)$. Feasibility flags and route counts can be cached because
they do not depend on the current LP solution. For a feasibility test alone,
cost and the visit-list tie-break may be omitted: resource dominance is
sufficient to preserve the existence of a feasible route. The stronger
frontier rule above is needed for the LA optimality argument.

The shortcutting assumption in Section~\ref{sec:vrptw-cuts} ensures that the
intersection of any feasible route with $S$ can itself be served by a
feasible route. For the VRPTW, nonnegative demands and travel times satisfying
the triangle inequality are sufficient, since arriving earlier can be
compensated by waiting. Thus the routes meeting $S$ induce a partition into
single-route-feasible subsets, and their number is at least $k(S)$.
Each such route exits every customer set containing $S$, establishing the
minimum-route inequality.

For a fixed support $\widehat C$, assign capacity $x_{ij}$ to each routing
edge, place both depot copies on the sink side, and add an artificial source
connected to each customer of $\widehat C$ by an infinite-capacity edge.
A minimum source--sink cut then returns the customer superset of
$\widehat C$ with minimum outgoing flow. Its inequality is separated when
this flow is below $k(\widehat C)$ by more than the separation tolerance.

\subsection{LA Separation and Cut Extraction}
\label{app:lm-separation}

Fix a support $C'\subseteq C$ and $\kappa=k(C')$. Let $\mathcal P$ be the
partitions of $C'$ into exactly $\kappa$ nonempty blocks. For a block $B$
with at least two customers, define
\[
 \mathcal E(B)=\bigcup_{u\ne v\in B}\mathcal E(u,B\setminus\{u,v\},v).
\]
For a singleton $B=\{u\}$, include the segment serving $u$, with no internal
edge and with both endpoints at $u$. Let
$\mathcal L=\bigcup_{\emptyset\ne B\subseteq C'}\mathcal E(B)$.
For $\gamma\in\mathcal L$, let $\xi_{ij}^\gamma$ indicate whether $\gamma$
uses internal edge $(i,j)$, and let $\xi_{i+}^\gamma$ and $\xi_{-i}^\gamma$
indicate whether it ends and starts at $i$, respectively. For the current
routing solution $x$, define
\begin{align*}
 d_{i+}(x)&=\sum_{j\notin C'}x_{ij},\\
 d_{-i}(x)&=\sum_{j\notin C'}x_{ji},\\
 b(x)&=\sum_{i\in C'}d_{i+}(x).
\end{align*}
Sums include only edges in $E$. The separation LP uses segment weights
$\Gamma_\gamma\ge0$, partition weights $\zeta_P\ge0$, and slack $s\ge0$:
\begin{subequations}
\label{eq:lm-la-primal}
\begin{align}
 \min_{\Gamma,\zeta,s\ge0}\quad &s\\
 \text{s.t.}\quad
 &\sum_{\gamma\in\mathcal L}\xi_{ij}^\gamma\Gamma_\gamma\le x_{ij}
 &&\forall (i,j)\in E\cap(C'\times C'),\label{eq:lm-la-internal}\\
 &\sum_{\gamma\in\mathcal L}\xi_{i+}^\gamma\Gamma_\gamma\le d_{i+}(x)
 &&\forall i\in C',\label{eq:lm-la-exit}\\
 &\sum_{\gamma\in\mathcal L}\xi_{-i}^\gamma\Gamma_\gamma\le d_{-i}(x)
 &&\forall i\in C',\label{eq:lm-la-entry}\\
 &\sum_{\gamma\in\mathcal E(B)}\Gamma_\gamma
 \ge\sum_{P\in\mathcal P:B\in P}\zeta_P
 &&\forall\emptyset\ne B\subseteq C',\label{eq:lm-la-block}\\
 &\sum_{P\in\mathcal P}\zeta_P+s\ge\kappa+1-b(x).
 \label{eq:lm-la-deficit}
\end{align}
\end{subequations}
The first three constraints limit segment weights by the available internal
and boundary flows. Each positive partition weight requires corresponding
segment weights for all its blocks. The last constraint requires sufficient
partition weight to cover the boundary deficit. The LP is always feasible,
and its optimal slack is zero when $b(x)\ge\kappa+1$. For $\kappa=1$, there
is only one partition; eliminating its weight leaves an LP in $\Gamma$ and $s$.

Introduce nonnegative dual variables $p_{ij}$, $p_{i+}$, $p_{-i}$ for
the flow bounds, $\mu_B$ for the block constraints, and $\rho$ for the
deficit constraint. The dual is given below, with its first constraint
imposed for every nonempty $B\subseteq C'$ and $\gamma\in\mathcal E(B)$:
\begin{align}
 \max\quad &\rho(\kappa+1-b(x))-\sum_{i,j\in C'}p_{ij}x_{ij}\nonumber\\
 &\quad-\sum_{i\in C'}\bigl(p_{i+}d_{i+}(x)+p_{-i}d_{-i}(x)\bigr)
 \label{eq:lm-la-dual}\\
 \text{s.t.}\quad
 &\mu_B\le\sum_{i,j\in C'}p_{ij}\xi_{ij}^\gamma\nonumber\\
 &\quad+\sum_{i\in C'}\bigl(p_{i+}\xi_{i+}^\gamma+p_{-i}\xi_{-i}^\gamma\bigr)
 \nonumber\\
 &\rho\le\sum_{B\in P}\mu_B\qquad\forall P\in\mathcal P,\nonumber\\
 &0\le\rho\le1,\qquad p,\mu\ge0.\nonumber
\end{align}
If the optimal slack is positive, an optimal dual solution yields
\begin{equation}
 \begin{aligned}
 &\sum_{i,j\in C'}p_{ij}x_{ij}+\rho b(x)\\
 &\quad+\sum_{i\in C'}\bigl(p_{i+}d_{i+}(x)+p_{-i}d_{-i}(x)\bigr)\\
 &\hspace{4em}\ge(\kappa+1)\rho.
 \end{aligned}
 \label{eq:lm-la-cut}
\end{equation}
The dual feasible region does not depend on $x$. Hence weak duality makes
the right-hand side minus the left-hand side of~\eqref{eq:lm-la-cut} a lower
bound on the optimal separation slack for every $x$. The inequality remains
applicable after the routing solution changes.
Using the routing degree equations, the boundary terms can equivalently be
eliminated by $d_{i+}(x)=1-\sum_{j\in C'}x_{ij}$ and
$d_{-i}(x)=1-\sum_{j\in C'}x_{ji}$.

\subsection{Preservation of a Common Optimum}

\begin{proposition}
\label{prop:lm-la-preservation}
Under the shortcutting assumption and the frontier rule of
Appendix~\ref{app:lm-frontier}, all inequalities~\eqref{eq:lm-la-cut}
preserve a common optimal routing solution.
\end{proposition}
\begin{proof}{Proof}
Start with an optimal routing solution. For each route, select the
lexicographically least pair of cost and visit list among feasible routes
serving the same customer set. These replacements preserve feasibility and
cannot increase cost, so the resulting solution $\hat x$ is optimal. Each
selected route belongs to its frontier, and Lemma~\ref{lem:lm-frontier}
places every contiguous subsegment in its corresponding frontier.

Fix any support $C'$. Its maximal contiguous runs in $\hat x$ partition
$C'$ and number $b(\hat x)\ge\kappa$. If $b(\hat x)\ge\kappa+1$, zero
segment and partition weights give zero slack in~\eqref{eq:lm-la-primal}.
Otherwise $b(\hat x)=\kappa$. Give weight one to the partition formed by
the runs and to the frontier segment corresponding to each run, using the
singleton convention when needed. These weights satisfy all constraints
with $s=0$. Weak duality therefore implies that $\hat x$ satisfies every
inequality~\eqref{eq:lm-la-cut}. The same $\hat x$ works for every support,
so the complete family retains the original optimal value. \Halmos
\end{proof}

\subsection{Tie-Breaking and the Reported Implementation}
\label{app:lm-frontier-audit}
The ordered visit list $\omega(\sigma)$ in
Appendix~\ref{app:lm-frontier} is compared lexicographically using fixed
customer indices. The first unequal index decides the order. Two segments
of the same class have the same endpoints and customer set, so adding a
common prefix or suffix preserves this order. This is the role of the
``word'' convention: it is simply the sequence of customer indices.

Resource dominance alone does not provide the strict comparison required
by Lemma~\ref{lem:lm-frontier}. For example, consider unit travel times and
costs, zero service times, upper time bounds of 100, and lower time bounds
of zero except $e_1=10$ and $e_4=20$. The two segments $(0,1,2,3)$ and
$(0,2,1,3)$ have equal cost and consumption, but earliest completion times
12 and 11, respectively. Their latest feasible starting times are both 97.
Resource dominance would discard the first. After appending $(3,4)$, both
completion times are 20 and both latest starting times are 96: the strict
resource improvement has disappeared. The first extended segment is
lexicographically smaller, although its prefix was discarded. All these
segments are feasible. This example explains why the resource-only rule
does not justify closure under taking subsegments.

The original implementation pruned an equal-cost segment on a strict
resource improvement, and used untruncated distance only when all compared
label components agreed. This differs from the rule used in
Proposition~\ref{prop:lm-la-preservation}. The implementation audit reported
33,712 missing trailing clips among 2,162,496 class checks in six
50-customer runs, and 1,814 classes with no retained segment having both
clips present. These checks establish a failure of the closure property,
not an observed invalid bound. Across the 168 original runs, no reported
lower bound exceeded the reference optimum and all 146 proved solutions
matched that optimum. These observations do not prove preservation of an
optimum in every instance.

A separate paired comparison of the two rules retained 158 instances on
which both root procedures converged and excluded ten that reached the
root time limit. The median root-gap change was zero and 75 pairs had
identical bounds. The largest reported improvement was RC204 with
25 customers, whose gap decreased from 8.66\% to 6.31\%; reported
worsenings were at most 0.20 percentage points. With a 60-second final-MILP
budget, both rules proved 125 instances, neither proved 30, and only the
corrected rule proved the remaining three. These are separate diagnostic
runs, not replacements for the 32-thread experiments in
Section~\ref{sec:experiments}. They indicate limited aggregate sensitivity,
but do not establish that the two rules are equivalent. The proof uses the
visit-order rule; the benchmark tables retain the original implementation.

\subsection{Separation Controls}
\label{app:lm-cut-controls}

\paragraph{Support-set generation.}
The separators consider customer subsets drawn from local clusters built
from the current fractional solution $x$. For each customer $u$, start
with $C_u=\{u\}$ and repeatedly add a customer $v\notin C_u$ maximizing
\[
 \sum_{w\in C_u}\phi_{wv}(x),
 \qquad
 \phi_{ij}(x)=
 \begin{cases}
  0,&x_{ij}\le\theta,\\
  x_{ij},&\theta<x_{ij}\le\tfrac12,\\
  1-x_{ij}+\nu',&x_{ij}>\tfrac12.
 \end{cases}
\]
Stop when $|C_u|=\bar m$ or the largest available gain is below $\nu$.
The score favors fractional connections, while $\nu'=1.01\nu$ allows a
single integral arc to extend a cluster. The implementation uses
$\bar m=12$, $\theta=0.01$, and $\nu=0.05$.

The nonempty subsets of these clusters supply the candidate supports for
both cut families, subject to the separator-specific size limits.
Repeated subsets need only be considered once. This construction targets
sets with substantial internal flow, for which the boundary flow may be
too small to accommodate the required routes or route segments. It is a
heuristic selection of supports: a larger cluster exposes more candidates,
at the cost of up to $2^{|C_u|}-1$ subsets per cluster. Restricting the
candidate pool affects which violated inequalities are found, not the
validity or common-optimum property established for each cut family.

\paragraph{Screening and limits.}
Candidate supports and their frontiers can be retained across iterations.
The boundary test $\kappa+1-b(x)>\epsilon_{\rm sep}$ avoids solving a
separation LP that cannot produce a cut above the violation tolerance.
Support-size and per-round cut limits restrict the candidates considered;
they do not affect the guarantee for cuts that are generated.

A frontier-size cap must be applied by skipping any support whose required
frontier exceeds the cap. Truncating a frontier and separating over the
remaining segments would not retain the witness used in
Proposition~\ref{prop:lm-la-preservation}. If route counts are capped at a
value $\bar k$, the resulting quantity $\min\{k(S),\bar k\}$ is a lower
bound rather than the exact count. It remains valid on the right-hand side
of a minimum-route inequality. LA separation can likewise use a certified
integer lower bound in place of $\kappa$: if an integer solution has more
than that many exits, its boundary deficit is already nonpositive.

The dual may admit several optima. Among optimal dual solutions, minimizing
$\sum_{i,j}p_{ij}+\sum_i(p_{i+}+p_{-i})$ gives a precise secondary
selection criterion. Equivalently, this criterion is induced
lexicographically by increasing all flow-bound right-hand sides by a common
positive perturbation. With a finite perturbation, the extracted row should
be checked for violation at the unperturbed $x$ before it is added.
The cut follows from dual feasibility; selection among dual solutions affects
its strength, not the common-optimum argument.

\clearpage
\section{Supplementary Cut Ablation}
\label{sec:exp-cuts}\label{app:lm-cut-ablation}

We assess the contribution of cuts through a separate, paired experiment
with one solver thread. The two configurations differ only in whether
root separation is enabled; aggregation, refinement, and reduced-cost removal
remain active in both. The experiment therefore measures the combined effect
of the cut families within DGAR, including their interaction with refinement
and pruning. Both configurations complete the root procedure on 166 instances.
The cut families are enabled or disabled together; this experiment does not
isolate their individual contributions.

On the 166 completed pairs, enabling cuts lowers the mean gap from 1.86\%
to 0.83\%, a 55.43\% reduction. Ninety-seven bounds improve beyond the
comparison tolerance and 69 remain within it; none is weaker
beyond this tolerance. The largest gain is on RC, whose mean falls from
4.21\% to 1.26\% (Table~\ref{tab:exp-ablation}). For RC1 at 50 customers,
the mean falls from 8.99\% to 0.09\%. Separation thus contributes most
strongly on the same subfamily for which DGAR performs particularly well in
the main comparison. The measured reductions apply to this one-thread ablation.

\begin{table}[!htbp]
\centering
\caption{Matched DGAR cut ablation at one thread, on 166 completed pairs.
Reduction is the relative reduction of the family mean. Improved/tied uses
the comparison tolerance.}
\label{tab:exp-ablation}\label{tab:lm-cut-ablation}
\small
\begin{tabular}{lrrrrr}
\toprule
Family & $n$ & No cuts gap (\%) & Cuts gap (\%) & Reduction (\%) & Improved/tied \\
\midrule
C & 51 & 0.04 & 0.02 & 31.38 & 5/46 \\
R & 69 & 1.64 & 1.14 & 30.91 & 56/13 \\
RC & 46 & 4.21 & 1.26 & 70.02 & 36/10 \\
\midrule
All & 166 & 1.86 & 0.83 & 55.43 & 97/69 \\

\bottomrule
\end{tabular}
\end{table}

Cuts can also change how much of the graph is needed. In
Figure~\ref{fig:exp-cut-scatter}, 95 pairs have fewer final arc-flow variables
with cuts, 61 have more, and ten are unchanged. The overall median paired
count ratio is 0.97. On RC1 with 50 customers it is 0.11: strong bounds
coincide with much smaller retained graphs. Cuts therefore have two potential
benefits: improving the relaxation and reducing the model passed to integer
optimization. The size effect depends on the instance, since cuts influence
both the pruning decisions and the refinement trajectory.

\begin{figure}[!htbp]
\centering
\includegraphics[width=\linewidth]{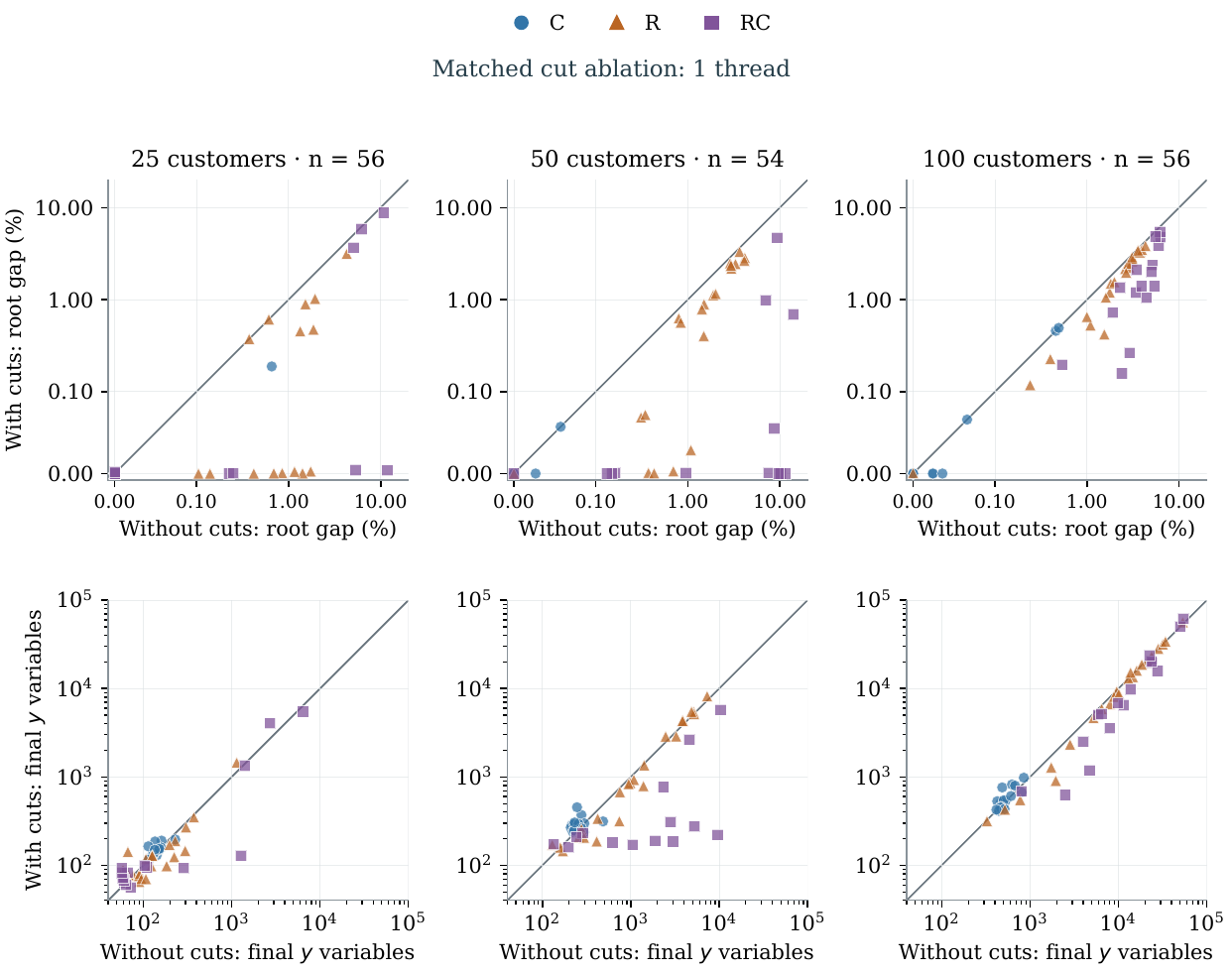}
\caption{DGAR cut ablation: one thread, 166 completed pairs. Horizontal axes
are without cuts and vertical axes with cuts. Top: root gap; bottom: final
arc-flow variable count. Points below the diagonal have a stronger bound or
a smaller graph with cuts, respectively. Shapes distinguish C, R, and RC.
Gap axes are linear up to 0.10\% and logarithmic thereafter; count axes
are logarithmic.}
\label{fig:exp-cut-scatter}
\end{figure}

\FloatBarrier

\clearpage
\section{Experimental Coverage and Measurement Details}
\label{app:exp-audit}

\subsection{Coverage and Run Selection}

The main 32-thread DGAR, compact, and RouteOpt root comparisons contain all
168 instances. DGAR uses the \texttt{pass2} configuration from source revision
\texttt{cc89eaa}. All 168 root procedures terminate with the recorded reason
\texttt{no\_progress}; none stops at the root time limit. This implementation
stopping condition does not independently establish equality with the fully
disaggregated LP bound, which was not computed for these instances. Bound comparisons use a
tolerance of 0.05 objective units to accommodate RouteOpt's one-decimal
output. Negative gaps within numerical tolerance are displayed as zero.
All calculations use unrounded values; continuous results are reported to
two decimal places.

BP-NG5 has 131 valid completed roots and BP-NG9 has 130. Their common valid
set contains 128 instances. BP-NG5 has 37 unfinished roots; BP-NG9 has 35
unfinished roots and three completed objectives above the reference value
(RC207\_50, R208\_100, and RC208\_100). These observations are excluded from
successful-root comparisons and listed in the accompanying per-instance
files. A final MIP bound is never substituted for a missing root bound.
The BP-NG runs used the same workstation as the other methods. Their solver
thread setting is listed as the Gurobi default in the treatment index and is
not reported in the per-instance outputs. They are used for the bound
and integer comparisons; the matched 32-thread root-time comparison includes
DGAR, compact, and RouteOpt.

The cut ablation uses the matched one-thread arms A and C from revision
\texttt{f74e8a2}. Their non-cut settings agree, while the four switches
controlling capacity and local-area separation are disabled in C. The
one-thread configuration also differs from the main 32-thread configuration
in LP restriction and LP algorithm settings; it is not a thread-only scaling
experiment. Two uncut runs reach their iteration limit and are excluded from
the 166 completed pairs: RC204\_50 has a partial uncut gap of 8.14\% and
RC208\_50 a gap of 13.96\%. Their completed cut-enabled one-thread gaps are
4.30\% and 9.24\%, respectively. Both instances remain in the complete
32-thread comparison.

The one-thread cut records contain 19 small downward lower-bound changes,
with maximum magnitude 0.009 objective units. These are within the stated
0.05-unit comparison tolerance; the figures and data retain the original
values rather than enforcing monotonicity.

\subsection{Timing and Model-Size Definitions}
\label{app:exp-time}

DGAR root wall time is the recorded total program time less the final-MILP
stage timer. The latter includes the final model construction and solve.
The root measure thus includes initial setup, input graph construction,
and the compression loop. Main/projector LP-call wall time is the sum of
\texttt{lower\_bound\_lp\_solve\_wall} and
\texttt{projector\_lp\_solve\_wall}. The local-area category combines the
LA-NG and three-customer separation timers, counting the aliased timer only
once. Other root work is the exact residual. Solver runtime counters are
reported separately and do not capture every separator solve; overlapping
stage timers are never added together.

The final $y$ count is \texttt{final\_milp\_edge\_flows}, separately from
\texttt{final\_milp\_actions} and the total number of model columns. Counts
are taken after closing cleanup and before the integer solver's presolve.
The last bound-producing LP column count comes from its own iteration record.
No earlier $y$ count is inferred by subtracting final routing variables from
an earlier total column count.
At 25, 50, and 100 customers, the median final row counts are 388.50,
1,004.00, and 5,450.00, respectively, over all 56 instances of each size.

The four-stage size comparison uses total model columns throughout. The
original count is the \texttt{root\_lp} entry in \texttt{REST\_SUMMARY}; this
model is constructed but not optimized. The first and final lower-bound LP
counts and objectives come from the first and last \texttt{compression}
records in \texttt{iteration\_records}, excluding projection-LP records.
The final MILP count comes from \texttt{final\_milp\_sizes}. These are
model snapshots before Gurobi presolve. The original networks include
artificial depot connectors; their removal is part of the final cleanup.
The first- and final-LP gaps use the same reference objectives and all 168
instances. The accompanying per-instance file retains all four size
snapshots, the two LP bounds, and their source hashes.

For integer outcomes, a proof requires the recorded solver optimal status
and agreement with the reference objective within 0.05 cost units. Table~\ref{tab:exp-integer}
and Figure~\ref{fig:exp-integer-profile} additionally require completion within
3,600 seconds. DGAR uses total program time; BP-NG5 and BP-NG9 use their overall
run timer; RouteOpt uses its reported elapsed time; Compact uses its MILP
optimizer timer. These are completion times, not first-incumbent times.
The comparison includes the complete RouteOpt search specified in
Section~\ref{sec:exp-setup}.

Four BP-NG5 and eight BP-NG9 optimal-status records fail the reference-objective
check, including one BP-NG5 record without a finite incumbent. They are not
counted as proofs. One otherwise qualifying BP-NG5 run and four BP-NG9 runs
finish after 3,600 seconds and are excluded from the common reporting window.
Their completion times range from 3,601.01 to 3,635.77 seconds. End-of-run
counts, before this time filter, are 88 and 96, respectively. All records and
exclusion reasons are retained in the accompanying data. Feasible objectives
remain available for time-limited runs.
DGAR's reference upper bound guides pruning and is distinct from its
feasible backup route. The backup matches the reference on 103 instances
and has higher cost on 65; the reference value is supplied information.

\subsection{Provenance and Reproduction}

The data come from the October 8, 2026 convergence archive, at repository
revision \texttt{33c550e}. The RouteOpt LP-thread rerun was performed on
October 3--4; the compact baseline was run with one worker and 32 solver
threads. The accompanying run records specify the configurations and solver settings.

RouteOpt uses version 2.0 with 32 LP threads and the configuration described
in Section~\ref{sec:exp-setup}. The executable checksum and full solver logs
are included in the accompanying data.

The accompanying package retains the raw JSON, comparator logs, source
hashes, common reference objectives, explicit common-set flags, and all
scripts needed to rebuild the tables and figures. The analysis recomputes
gaps from objectives rather than reading rounded spreadsheet gaps. The
one-thread cut data are stored separately from the 32-thread main experiment.
\FloatBarrier

\subsection{Paired Root Comparisons}
\label{app:exp-scatter}

Figure~\ref{fig:exp-paired-root} compares the root bounds and elapsed root
times for DGAR and RouteOpt on all 168 instances. DGAR occupies the vertical
axis. Points below the diagonal indicate a smaller gap or shorter root time,
respectively. The size and family breakdown distinguishes the bound gains
from the computation needed to obtain them.

\begin{figure}[!htbp]
\centering
\includegraphics[width=\linewidth]{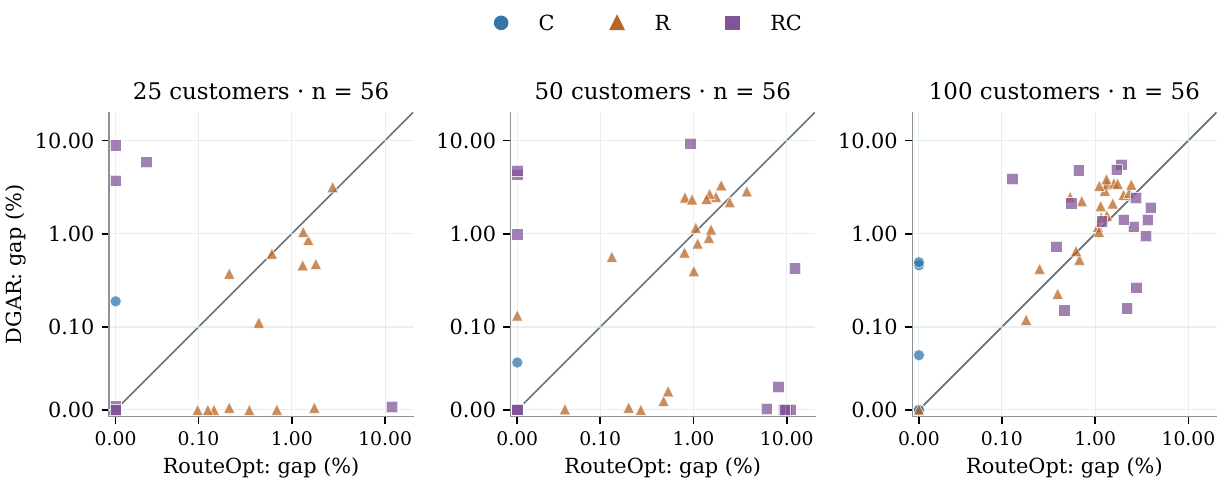}
\smallskip
\includegraphics[width=\linewidth]{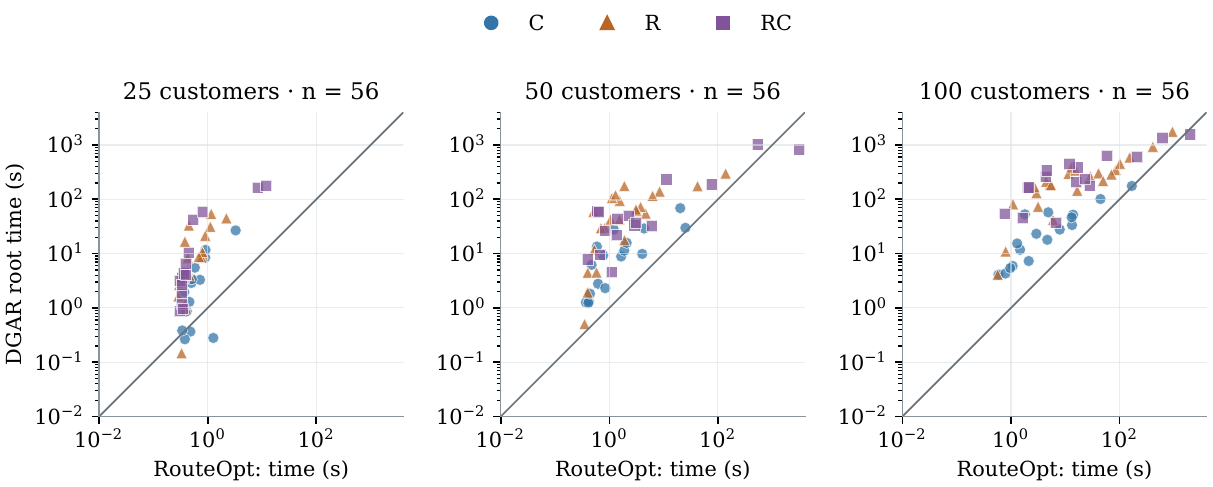}
\caption{Paired root comparisons at 32 threads. Top: root gaps; bottom:
elapsed root times. All 56 instances are included at each size. Symbols
distinguish C, R, and RC. The gap axes are linear up to 0.10\% and logarithmic
thereafter; the time axes are logarithmic. Timing definitions are given in
Appendix~\ref{app:exp-time}.}
\label{fig:exp-paired-root}
\end{figure}
\FloatBarrier

}
\end{APPENDICES}

\end{document}